\documentclass[12pt]{article}
\usepackage{amsmath,amssymb,amsthm,amsfonts,graphicx,authblk,xcolor}
\begin{document}

  \newtheorem{theorem}{Theorem}[section]
  \newtheorem{proposition}[theorem]{Proposition}
  \newtheorem{lemma}[theorem]{Lemma}
  \newtheorem{definition}[theorem]{Definition}
  \newtheorem{cor}[theorem]{Corollary}
  \newtheorem{remark}[theorem]{Remark}
  \newtheorem{algorithm}[theorem]{Algorithm}
\newcommand{\rdots}{.\hspace{.1em}\raisebox{.8ex}
                    {.\hspace{.1em}\raisebox{.8ex}{.}}}
\newcommand{\sC}
           {\mbox{\rm C\hspace{-0.4em}\rule{0.05ex}{1.55ex}\hspace{0.4em}}}
\newcommand{\sN}
           {\mbox{\rm N\hspace{-0.71em}\rule{0.05ex}{1.55ex}\hspace{0.71em}}}
\newcommand{\sR}
           {\mbox{\rm R\hspace{-0.7em}\rule{0.08ex}{1.55ex}\hspace{0.7em}}}
\newcommand{\sT}
           {\mbox{\rm T\hspace{-0.5em}\rule{0.05ex}{1.55ex}\hspace{0.5em}}}

\renewcommand{\theequation}{\arabic{section}.\arabic{equation}}

\newcommand{\be}{\begin{equation}}
\newcommand{\ee}{\end{equation}}
\newcommand{\ba}{\begin{array}}
\newcommand{\ea}{\end{array}}
\newcommand{\C}{{\mathbb C}}
\newcommand{\cH}{{\cal H}}
\newcommand{\al}{{\alpha}}
\newcommand{\bt}{{\beta}}
\newcommand{\G}{{\Gamma}}
\newcommand{\g}{{\gamma}}
\newcommand{\dl}{{\delta}}
\newcommand{\Dl}{{\Delta}}
\newcommand{\co}{{\rm const}\,}
\newcommand{\comment}[1]{}
\newcommand{\ds}{\displaystyle}
\newcommand{\eop}{\quad\rule{7pt}{8pt}\vspace*{3mm}}
\newcommand{\eq}[2]{\begin{equation}\label{#1}#2\end{equation}}
\newcommand{\ha}{\frac{1}{2}}
\newcommand{\La}{\Lambda}
\newcommand{\la}{\lambda}
\newcommand{\vp}{\varphi}
\newcommand{\tl}{\tilde}
\newcommand{\ma}[2]{\left[\begin{array}{#1}#2\end{array}\right]}
\newcommand{\Th}{\Theta}
\newcommand{\om}{\omega}
\newcommand{\Om}{\Omega}
\newcommand{\sg}{\sigma}
\newcommand{\sq}{\sqrt}
\newcommand{\txt}[1]{\quad{\rm #1}\quad}
\newcommand{\wh}{\widehat}
\newcommand{\wt}{\widetilde}

\newcommand{\canc}[1]{\textcolor[rgb]{.7,.7,.7}{#1}}
\newcommand{\ins}[1]{\textcolor[rgb]{1,0,0}{#1}}
\newcommand{\revpaola}[1]{\textcolor{violet}{#1}}

\title{Quasiseparable LU decay bounds for inverses of banded matrices}

\date{}
\author[1]{P. Boito} 
\author[2]{Y. Eidelman}
\affil[1]{Dipartimento di Matematica, Universit\`a di Pisa, Largo Bruno Pontecorvo, 5 - 56127 Pisa, Italy. Email: paola.boito@unipi.it }
\affil[2]{School of Mathematical Sciences, Raymond and Beverly
Sackler Faculty of Exact Sciences, Tel-Aviv University, Ramat-Aviv,
69978, Israel. Email: eideyu@tauex.tau.ac.il}

\maketitle
\thispagestyle{empty}
\abstract{We develop new, easily computable exponential decay bounds for inverses of banded matrices, based on the quasiseparable representation of Green matrices. The bounds rely on a diagonal dominance hypothesis and do not require explicit spectral information. Numerical experiments and comparisons show that these new bounds can be advantageous especially for nonsymmetric or symmetric indefinite matrices.}

\section{Introduction}
A well-known result in matrix theory states that, under suitable hypotheses, the inverse of an $n\times n$ banded matrix $A$, while not being banded itself in general, exhibits an exponential off-diagonal decay behavior. In other words, there exist constants $K>0$ and $0<\xi<1$, independent of $n$, such that 
\begin{equation}\label{generalbound}
|A^{-1}(i,j)|\leq K \xi^{|i-j|}.
\end{equation}
It is worth emphasizing that, although this result is often stated for a single matrix $A$, what we really have in mind here is a sequence of invertible matrices $\{A_n\}$ of increasing size and uniform bandwidth: then the bound \eqref{generalbound} holds uniformly for all inverses $\{A_n^{-1}\}$, precisely because $K$ and $\xi$ do not depend on $n$. Clearly, an {\em a priori} bound like \eqref{generalbound} is especially useful when both $K$ and $\xi$ are explicitly (and easily) computable. Existence results (see, e.g., \cite{Groechenig06, Jaffard}) are of remarkable theoretical interest, but are obviously difficult to apply in practice.
 
Starting from Demko, Moss and Smith's seminal paper \cite{DMS84}, decay results in the flavor of \eqref{generalbound} have been proved in many versions. The original bound in \cite{DMS84} relies on polynomial approximation of the function $x\rightarrow 1/x$ and holds for Hermitian, positive definite matrices whose spectrum is contained in a fixed interval $[a,b]$, with $0<a<b$. The constants $K$ and $\xi$ essentially depend on the condition number $b/a$ and on the bandwidth of the matrix. 

This bound has been subsequently extended to more general functions and sparsity patterns \cite{Baskakov, BenziRazouk07, Mastronardi, BenziBoito14, BenziSimoncini, FrommerNLAA, Frommer18, Nabben, Pozza}. A recent contribution concerning spectral projectors of Hermitian matrices is found in \cite{BenziRinelli}. Note that all these results require spectral information on $A$. This can sometimes be obtained cheaply through Gershgorin's theorems, but loose estimates on the eigenvalues of $A$ typically yield pessimistic decay bounds.

Thanks to the spectral theorem, approximation-based bounds for $f(A)$ are much easier to obtain when $A$ is Hermitian, or, more generally, normal. A straightforward generalization to non-normal, diagonalizable $A$ is possible, but it involves the condition number of the eigenvector matrix, which can be large and difficult to compute. Another approach relies on approximation over the field of values of $A$; see \cite{Crouzeix} and references therein. In fact, a common trait of most approximation-based bounds for $f(A)$ is the requirement that $f(x)$ must be analytic on a suitable convex subset of the complex plane containing the spectrum of $A$, for instance a Bernstein ellipse or the field of values of $A$. 

As a consequence, such bounds 
cannot be applied directly to inverses of matrices containing both positive and negative real eigenvalues. One exception is the work by Chui and Hasson \cite{CH83, Hasson07}, which yields however a decay rate but not computable decay bounds. A possible workaround consists in rewriting $A^{-1}$ as $A^H(A A^H)^{-1}$, as in \cite{DMS84}, but this approach considerably worsens decay rates and cannot be applied to one-sided banded matrices, since in this case the banded structure would be lost. 

In this work we revisit decay bounds for inverses of banded matrices from a different viewpoint. Banded matrices and their inverses belong to the class of quasiseparable matrices, for which a compact representation with $O(n)$ parameters is available \cite{EGH1}. It is therefore natural to ask: given a banded matrix $A$, can we exploit a quasiseparable representation of $A^{-1}$ to formulate {\em a priori} decay bounds such as \eqref{generalbound}? 
Similar research directions have been abundantly explored for tridiagonal $A$, see e.g., \cite{Dedieu88, Frommer18, Meurant, Nabben},
but little work is available for larger bandwidths \cite{Nabben}. 

Here we rely on the results presented in \cite{BE23}, which provide quasiseparable representations of Green matrices and, in particular, of inverses of one-sided and two-sided banded matrices. The hypotheses on the properties and location of the spectrum of $A$ are replaced by a strong dominance condition on $A$, which can be easily checked. The representation of $A^{-1}$ is derived from a structured LU factorization of $A$ and the decay bounds on $A^{-1}$ are obtained through a careful analysis of the magnitude of elements in the Green generators. We note in passing that this approach also involves a characterization of off-diagonal decay in the factors of the LU factorization of $A$, which is a research direction of interest in its own right, for instance for preconditioning purposes \cite{Bellavia, Krishtal}.

The new bounds presented in this work have the merit of being easily computable, without the need for explicit spectral estimates for $A$. Numerical experiments show that in many cases they can give a better description of the decay behavior of $A^{-1}$ than other bounds available in the literature, particularly when $A$ is symmetric indefinite or nonsymmetric.  
 
The paper is organized as follows. Section \ref{sec:background} recalls definitions, notation and the relevant results from \cite{BE23}. The main result on decay bounds, namely, Theorem \ref{TTR}, is presented in Section \ref{sec:main}, whereas Section \ref{sec:proof} is devoted to the proof of Theorem \ref{TTR}, together with auxiliary results and comments. Section \ref{sec:QR} briefly touches on the development of decay bounds from a structured QR factorization of $A$. Numerical experiments and comparisons with bounds from the literature are presented in Section \ref{sec:numerical}, whereas Section \ref{sec:conclusions} contains concluding remarks and ideas for further investigation.

% Say somewhere (abstract, intro, cover letter...) that these results are new and can be of interest for two main reasons: they provide a new approach to the development of decay bounds, which had previously been explored essentially only for tridiagonal matrices, and the resulting bounds are more effective than previous ones on many examples, particularly on nonsymmetric or on symmetric indefinite matrices.

% What about operators?

\section{Background}\label{sec:background}
We recall here the main definitions and results on quasiseparable and Green matrices that will be used later. Details and proofs can be found in \cite{BE23} and \cite{EGH1}. Note that throughout this paper we use MATLAB notation to denote submatrices, that is, $A(i:j,k:\ell)$ denotes the submatrix of $A$ defined by row indices from $i$ to $j$ and column indices from $k$ to $\ell$.

\begin{definition}
Let $r, N$ be integers such that $N>r>0$.

An $N\times N$ scalar matrix $A=\{A(i,j)\}_{i,j=1}^N$ is called {\em a lower
band matrix of order $r$}  if $A(i,j)=0$ for $i-j>r$.

A matrix $B$ is called {\em a lower Green matrix of order $r$} if
\be\label{natu}
{\rm rank}\,B(k:N,1:k+r-1)\le r,\quad k=1,2,\dots,N-r.
\end{equation}

A matrix $D$ is called {\em an upper Green matrix of order $r$} if
\be\label{natuu}
{\rm rank}\,D(1:k+r-1,k:N)\le r,\quad k=1,2,\dots,N-r.
\end{equation}
\end{definition}

Banded and Green matrices are deeply connected: a well-known result by Asplund \cite{Asplund59} states that the class of invertible lower Green matrices of order $r$
coincides with the class of inverses of invertible lower band
matrices of the same order $r$.

In this work we will employ a particular quasiseparable representation for Green matrices, which we summarize below; see \cite{EGH1} and \cite{BE23} for details.
Let $F=\{F'(i,j)\}_{i,j=1}^K$ be a block matrix
with entries of sizes $m_i\times n_j$. Throughout this paper, we use the prime $'$ symbol to emphasize the fact that we are referring to a block partition of a matrix: for instance, $F(1,1)$ denotes the (scalar) entry of $F$ in position $(1,1)$, whereas $F'(1,1)$ denotes the block of $F$ in position $(1,1)$, which has size $m_1\times n_1$. 

Recall that, if $F$ has a lower quasiseparable structure (i.e., submatrices in the strictly lower triangular part of $F$ have low rank), the strictly lower
triangular part of $F$ admits the quasiseparable representation
\be\label{qrprl}
F'(i,j)=p(i)a_{ij}^{>} q(j),\;1\le j<i\le K,
\end{equation}
where $p(i),\;i=1,\dots,K,\;q(j),\;j=1,\dots,K-1,\;
a(k),\;k=2,\dots,K-1$ are matrices of (small) sizes $m_i\times r_{i-1},r_i\times n_j,
r_k\times r_{k-1}$, respectively. Here we denote $a_{ij}^{>}=a(i-1)\cdots a(j+1)$, with the convention that 
$a_{ij}^{>}$ is the identity matrix if $j\geq i-1$.

In the {\em Green representation} that will be used in this paper, we treat an $N\times N$ scalar matrix as a block matrix of size $(N-r+2)\times
(N-r+2)$. The sizes of the blocks are chosen as
\be\label{msh}
\begin{gathered}
m_0=0,\;m_1=m_2=\dots=m_{N-r}=1,\;m_{N-r+1}=r;\\
n_0=r;\;n_1=n_2=\dots=n_{N-r}=1,\;n_{N-r+1}=0.
\end{gathered}
\end{equation}
Note that indices here start from zero.
A lower band matrix of order $r$ can be seen as a block matrix
with $n_i\times m_j,\;i,j=0,1,\dots,N-r+1$ entries and therefore turns out to be block upper
triangular. A lower Green matrix of order $r$ is seen as a block one with
entries of sizes $m_i\times n_j,\;i,j=0,1,\dots,N-r+1$.

It can be shown that a lower Green matrix $B$ of order $r$ can be written in block form as
\begin{equation}\label{eq:greenrep}
B'(i,j)=p(i)a^>_{ij}q(j),\quad 0\le j<i\le N-r+1
\end{equation}
with matrices $p(i),\;i=1,\dots,N-r+1$ of sizes $m_i\times r$,
$q(j),\;j=0,\dots,N-r$ of sizes $r\times n_j$ and $a(k),\;k=1,\dots,N-r$
of sizes $r\times r$. Without loss of generality we
will assume $q(0)=I_r$.

Thus the block lower triangular part of a lower Green matrix of order $r$, that is, the part of the matrix where scalar indices $(i,j)$ are such that $j-i\le r$, is completely defined by the parameters
$p(i),q(i),a(i),\;i=1,\dots,N-r,\;p(N-r+1)$.
For instance an $(r+4)\times(r+4)$ matrix $B$ has the form
\begin{gather*}
B'=\\
\left(\ba{ccccc}p(1)&\ast&\ast&\ast&\ast\\
p(2)a(1)&p(2)q(1)&\ast&\ast&\ast\\
p(3)a(2)a(1)&p(3)a(2)q(1)&p(3)q(2)&\ast&\ast\\
p(4)\cdots a(1)&p(4)a(3)a(2)q(1)&p(4)a(3)q(2)&p(4)q(3)&\ast\\
p(5)\cdots a(1)&p(5)a(4)a(3)a(2)q(1)&p(5)a(4)a(3)q(2)&
p(5)a(4)q(3)&p(5)q(4)\ea\right).
\end{gather*}
In particular, note that the first column is a block column of width $r$ and the last row is a block row of height $r$, whereas the remaining entries are scalar.

The elements $p(i),q(i),a(i),\;i=1,\dots,N-r,\;p(N-r+1)$, where
$p(i)\;(i=1,\dots,N-r)$ are $r$-dimensional rows and $p(N-r+1)$ is an
$r\times r$ matrix, $q(i)$ are $r$-dimensional
columns and $a(i)$ are $r\times r$ matrices,
are said to be {\it lower Green generators} of the matrix $B$. Conversely, if the representation (\ref{eq:greenrep})
holds, then $B$ is a lower Green matrix of order $r$.

\begin{remark}\label{rem:green2scalar}
Let us write explicitly the conversion from the Green block representation \eqref{eq:greenrep} of $B$ to the usual scalar representation of the triangular part of $B$, which is indexed from $1$ to $N$. For $1\leq i\leq j\leq N$ we have
\begin{equation}\label{eq:green2scalar}
B(i,j)=\left\{\begin{array}{l}
\left[p(i)a_{i,0}^>\right](j)\qquad {\rm if }\; 1\leq i\leq N-r,\; 1\leq j\leq r,\\
p(i)a_{i,j-r}^> q(j-r) \qquad {\rm if }\; 1\leq i\leq N-r,\; r+1\leq j\leq N,\\
\left[p(N-r+1)a_{N-r+1,j-r}^>q(j-r)\right](i-N+r) \\ \hspace{2cm} {\rm if }\; N-r+1\leq i\leq N,\; r+1\leq j\leq N.
\end{array}\right.
\end{equation}
%In fact, in our Green representation we may also avoid the use of $r\times r$ generator $p(N-r+1)$ and introduce 
%row vectors $p(i)$ for $i=N-r+1,\ldots,N$, so that $p(N-r+\ell)$ coincides with row $\ell$ of matrix $p(N-r+1)$, for $\ell=1,\ldots,r$. If we adopt this modified definition (as we will do in the next section) we can drop the last line of \eqref{eq:green2scalar} and take $1\leq i\leq N$ in the second line. %Another option to encode the last block row is to define extra generators $p(i)$, $a(i)$, $i=N-r+1,\ldots,N-1$ of respective sizes $1\times r_i$, $r_i\times r_{i+1}$ where $r_i = N-i$.
\end{remark}

Let us now recall some results from \cite{BE23} that will be useful in the 
development of our decay bounds. The following theorem provides a structured LU 
factorization of a lower banded matrix.

\begin{theorem}\label{ICF1}
Let $A=\{A(i,j)\}_{i,j=1}^N$ be a strongly regular lower band matrix of order $r$.

Then $A$ admits the factorization
\be\label{irn15}
A=LR,
\end{equation}
where $L$ is a unit lower triangular matrix and $R$ is an upper triangular 
matrix. The inverse $L^{-1}$ of the lower triangular factor $L$ may be 
represented as the product
\be\label{un}
L^{-1}=\tilde L_{N-r+1}\tilde L_{N-r}\cdots\tilde L_2\tilde L_1
\end{equation}
with
\be\label{un1}
\tilde L_k=I_{k-1}\oplus L_k\oplus I_{N-k-r},\; k=1,\dots,N-r,\quad
\tilde L_{N-r+1}=I_{N-r}\oplus L_{N-r+1},
\end{equation}
and $(r+1)\times(r+1)$ lower triangular matrices
$L_k,\;k=1,\dots,N-r$ and $r\times r$ lower triangular matrix $L_{N-r+1}$.
Moreover, the lower triangular matrices $L_k$ as well as the upper triangular
entries of the matrix $R$ are obtained as follows.

1. Set
\be\label{lry0}
Y_0=A(1:r,1:N),\quad\gamma_1=Y_0(1,1).
\end{equation}

2. For $k=1,\dots,N-r$ perform the following.
Set 
\be\label{irks}
\gamma_k=Y_{k-1}(k,k),\; X_k=Y_{k-1}(k,k+1:N)
\end{equation}
and compute
\be\label{lnk}
f_k=\left(\ba{c}Y_{k-1}(k+1:k+r-1,k)\\A(k+r,k)\ea\right)\frac1{\gamma_k}.
\end{equation}
Set
\be\label{ol}
Z_k=\left[\ba{c}Y_{k-1}(k+1:k+r-1,k+1:N)\\A(k+r,k+1:N)\ea\right]
\end{equation}
and compute
\be\label{olz}
Y_k=-f_k\cdot X_k+Z_k.
\end{equation}
Set
\be\label{natas}
L_k=\left(\ba{cc}1&0\\-f_k&I_r\ea\right).
\end{equation}
Set
\be\label{olia}
R(k,k)=\gamma_k,\quad R(k,k+1:N)=X_k.
\end{equation}

3. For $k=N-r+1,\dots,N-1$ perform the following.

Set $\gamma_k=Y_{k-1}(k,k),\; X_k=Y_{k-1}(k,k+1:N)$ and compute
\be\label{lnk3}
f_k=Y_{k-1}(k+1:N,k+1)\frac1{\gamma_k}.
\end{equation}
Set
\be\label{ol3}
Z_k=Y_{k-1}(k+1:N,k+1:N)
\end{equation}
and compute
\be\label{olz3}
Y_k=-f_k\cdot X_k+Z_k.
\end{equation}
Set
$$
L_k=\left(\ba{cc}1&0\\-f_k&I_r\ea\right).
$$
Set
\be\label{olia3}
R(k,k)=\gamma_k,\quad R(k,k+1:N)=X_k.
\end{equation}

4. Set $\gamma_N=Y_{N-1},\;R(N,N)=\gamma_N$.
\end{theorem}

Quasiseparable generators for the lower triangular factor can be extracted from
 the matrices $L_k$ as follows.

\begin{lemma}\label{AA}
The matrix $L^{-1}$ in (\ref{un}), (\ref{un1}) is a lower Green - upper
band of order $r$ matrix with lower Green generators $p_L(i),q_L(i),a_L(i)\;
i=1,\dots,N-r$ and diagonal entries $d_L(k)\;(k=1,\dots,N-r)$ obtained from the
partitions
\be\label{iof5}
L^{-1}_k=\left[\ba{cc}p_L(k)&d_L(k)\\a_L(k)&q_L(k)\end{array}\right],\quad
k=1,\dots,N-r
\end{equation}
and by setting
\be\label{iof5r}
p_L(N-r+1)=L_{N-r+1}.
\end{equation}
\end{lemma}

Building on Theorem \ref{ICF1} and Lemma \ref{AA}, we can now compute 
quasiseparable generators for the lower triangular part of the inverse of a 
lower banded matrix.

\begin{theorem}\label{UF1}
Let $A$ be a strongly regular lower band matrix of order $r$.

Then lower Green generators $p(i),q(i),a(i)\;(i=1,\dots,N-r)$,
$P_{N-r+1}$ of the matrix $A^{-1}$ are obtained as follows.

1. Using the algorithm from Theorem \ref{ICF1} compute the lower triangular
 matrices
$L_k\;(k=1,\dots,$
$N-1)$ of orders $r_k=r$ for $k=1,\dots,N-r$ and
$r_k=N-k+1$ for $k=N-r+1,\dots,N-1$,
as well as diagonal entries and subrows
\be\label{tre}
\gamma_k=R(k,k),\;k=1,\dots,N,\quad X_k=R(k,k+1:N),\;k=1,\dots,N-1
\end{equation}
of the lower triangular matrix $R$.
Determine the lower Green generators $p_L(i),q_L(i)$,
$a_L(i)\;(i=1,\dots,N-r)$ of the 
lower Green-upper band matrix $L^{-1}$ via partitions (\ref{iof5}) and the
matrices $p_L(i),a_L(i)\;(i=N-r+1,\dots,N-1)$ of sizes $1\times r_i,
r_i\times r_{i+1}$ from the partitions
\be\label{iof5.}
L_i=\left[\ba{c}p_L(i)\\a_L(i)\end{array}\right],\quad
i=N-r+1,\dots,N-1.
\end{equation}

2. Compute the lower Green generators of the matrix $A^{-1}$ as follows.

2.1. Compute the lower Green generator $P_{N-r+1}$ as follows.
Set
\be\label{rn}
p(N)=P_N=\frac1{\gamma_N},
\end{equation}
and for $k=N-1,\dots N-r+1$ compute
\be\label{mishln}
p(k)=\frac1{\gamma_k}(p_L(k)-X_kP_{k+1}a_L(k))
\end{equation}
\be\label{mashln}
P_k=\left(\ba{c}p(k)\\P_{k+1}a_L(k)\ea\right).
\end{equation}

2.2. Compute the lower Green generators $p(k),q(k),a(k)$ as follows.

2.2.1. Set
\be\label{masl}
q(k)=q_L(k),\;a(k)=a_L(k),\;k=1,\dots,N-r.
\end{equation}

2.2.2. For $k=N-r,\dots 1$ compute
\be\label{mishl}
p(k)=\frac1{\gamma_k}(p_L(k)-X_kP_{k+1}a(k))
\end{equation}
\be\label{mashl}
P_k=\left(\ba{c}p(k)\\P_{k+1}a(k)\ea\right).
\end{equation}
\end{theorem}

\begin{remark}
The generator $P_{N-r+1}$ from step 2.1 of Theorem \ref{UF1} can also be computed as $P_{N-r+1}=R(N-r+1:N,N-r+1:N)^{-1}L_{N-r+1}$.
\end{remark}

\section{LU-based decay bounds: the main result}\label{sec:main}
\setcounter{equation}{0}
 
Let us now introduce the main result of this paper, that is, estimates for entries of the inverses of (possibly one-sided) band matrices in terms of the matrix elements. 
The main hypothesis here is a diagonal dominance condition, which is easy to check on a given matrix. 

\begin{theorem}\label{TTR}
Let $A=\{A(i,j)\}_{i,j=1}^N$ be an invertible lower-banded matrix of order $r$. 
Assume that $A$ satisfies the following strong dominance condition 
\begin{equation}\label{gl1}
\mu|A(k,k)|\geq \sum_{i=1}^{k-1}|A(i,k)|+\sum_{i=k+1}^{k+r}|A(i,k)|,\;
k=1,2,\dots,N
\end{equation}
with $0<\mu<1$.

Then for the Green block representation of $A^{-1}$ it holds
\begin{equation}\label{mrn}
|(A^{-1})'(i,j)|\le M\gamma^{i-j-r},\quad i - j > r
\end{equation}
with $\gamma=\mu^{1/r}$ and $M=\frac{1+\mu^2}{(1-\mu)(1-\mu^2)\min_i|A(i,i)|}.$
\end{theorem}

The proof of Theorem \ref{TTR} will be presented in the next section. 

Note that \eqref{mrn} is given in terms of the Green blocks of $A^{-1}$. For practical purposes it may be more useful to reframe this result in terms of the usual scalar representation. 

\begin{theorem}\label{thm:bound}
Let $A$, $\mu$, $\gamma$, $M$ be as in Theorem \ref{TTR}. Then for $i\geq j$ it holds
\begin{equation}\label{simplebound}
|(A^{-1})(i,j)|\le M\gamma^{i-j}.
\end{equation}
\end{theorem}

\begin{proof}
Consider at first the case $i=j$. Using Corollary \ref{cor:varah} in Section \ref{sec:proof} we have
$$
|(A^{-1})(i,i)|\leq\|A^{-1}\|_1<\frac{1}{(1-\mu)\min_k|A(k,k)|}<M.
$$
Let us now discuss the case $i>j$. Define the augmented matrix
$$
\mathcal{A}=I_r\oplus A\oplus I_r\in\mathbb{R}^{\tilde{N}\times\tilde{N}},
$$
with $\tilde{N}=N+2r$. Note that $\mathcal{A}$ satisfies hypothesis \eqref{gl1} with the same value of $\mu$ as $A$. Combining \eqref{mrn} applied to $\mathcal{A}$ with Remark \ref{rem:green2scalar} we can write
\begin{equation}\label{mrns}
|(\mathcal{A}^{-1})(i,j)|\le \left\{\begin{array}{l}
M\gamma^{i-r},\qquad r+1\leq i\leq \tilde{N}, 1\leq j\leq r\\
M\gamma^{i-j},\qquad r+1\leq j<i\leq \tilde{N}-r\\
M\gamma^{\tilde{N}-r+1-j},\qquad \tilde{N}-r+1\leq i\leq \tilde{N}, r+1\leq j \leq \tilde{N}-r.\\
\end{array}\right.
\end{equation}
The thesis follows from the second line in \eqref{mrns}.
\end{proof}

%This is done by combining \eqref{mrn} with Remark \ref{rem:green2scalar}.
%\begin{cor}
%Let $A$, $\mu$, $\gamma$, $M$ be as in Theorem \ref{TTR}. Then it holds
%\begin{equation}\label{mrns}
%|(A^{-1})(i,j)|\le \left\{\begin{array}{l}
%M\gamma^{i-r},\qquad r+1\leq i\leq N, 1\leq j\leq r\\
%M\gamma^{i-j},\qquad r+1\leq j<i\leq N.
%\end{array}\right.
%\end{equation}
%Moreover, \eqref{mrns} implies for $i-j>r$
%\begin{equation}\label{simplebound}
%|(A^{-1})(i,j)|\le M\gamma^{i-j-r}.
%\end{equation}
%\end{cor} 

\begin{remark} 
In the case of a two-sided banded matrix $A$ of order $r$, the 
condition (\ref{gl1}) takes the form
$$
\mu|A(k,k)|\geq \sum_{i=k-r}^{k-1}|A(i,k)|+\sum_{i=k+1}^{k+r}|A(i,k)|.
$$
\end{remark}

\section{Proof of LU-based decay bounds}\label{sec:proof} 
Let us start with two technical lemmas.

\begin{lemma}\label{lemma:fk}
Under the hypotheses of Theorem \ref{TTR}, it holds 
\be\label{larr}
\|f_k\|_1\le\mu<1,\quad k=1,\dots,N-r
\end{equation}
\end{lemma}

\begin{proof}
Consider the first step of the algorithm from Theorem \ref{ICF1}.
Note that the hypothesis (\ref{gl1}) implies
\be\label{simplegl1}
\mu|A(k,k)|\ge \sum_{i=1}^{k-1}|A(i,k)|+\sum_{i=k+1}^{k+r}|A(i,k)|,\;
k=1,2,\dots,N-r.
\end{equation}
This yields for $k=1$ 
$$
\mu|A(1,1)|\ge \|A(2:r+1,1)\|_1,
$$ 
which implies $\|f_1\|_1\le\mu$. As a result of the first step we obtain the 
matrix 
$$
{\mathcal A}_1=\tilde L_1A=\left(\ba{cc}\gamma_1&H_1\\0&A_1\ea\right),
$$
where $H_1=A(1,2:N)$ and the $(N-1)\times(N-1)$ matrix $A_1$ is the Schur complement
\be\label{ire1}
A_1=A(2:N,2:N)-f_1H_1.
\end{equation}
Next we apply the same procedure to $A_1$, which is lower banded of order $r$, and so on. In order to prove the lemma, it suffices to prove that the matrix $A_1=\{A_1(i,j)\}_{i,j=2}^N$ 
satisfies the diagonal dominance conditions (\ref{simplegl1}), i.e.
\be\label{iroc}
\mu|A_1(k,k)|\ge \sum_{i=2}^{k-1}|A_1(i,k)|+\sum_{i=k+1}^{k+r}|A_1(i,k)|,\;
k=2,\dots,N-r.
\end{equation}
Using (\ref{ire1}) we have
\be\label{qyr}
A_1(k,k)=A(k,k)-f_1(k)A(1,k),\;k=2,\dots,N-r
\end{equation}
and
\be\label{qyrj}
A_1(i,k)=A(i,k)-f_1(i)A(1,k),\;k=2,\dots,N-r,\;i=2,\dots,N-r,\;i\neq k.
\end{equation}
Using (\ref{qyr}) and (\ref{qyrj}) we obtain the inequalities
\be\label{gyr}
|A_1(k,k)|\ge|A(k,k)|-|f_1(k)||A(1,k)|,\;k=2,\dots,N-r
\end{equation}
and
\be\label{gyrj}
|A_1(i,k)|\le|A(i,k)|+|f_1(i)||A(1,k)|,\;
k=2,\dots,N-r,\;i=2,\dots,N-r,\;i\neq k.
\end{equation}
Summing the inequalities in (\ref{gyrj}) over index $i$ we get
\begin{gather*}
\sum_{i=2}^{k-1}|A_1(i,k)|+\sum_{i=k+1}^{k+r}|A_1(i,k)|\le\\
\sum_{i=2}^{k-1}|A(i,k)|+\sum_{i=k+1}^{k+r}|A(i,k)|+
\left(\sum_{i=2}^{k-1}|f_1(i)|+\sum_{i=k+1}^{k+r}|f_1(i)|\right)|A(1,k)|\le\\
\sum_{i=1}^{k-1}|A(i,k)|+\sum_{i=k+1}^{k+r}|A(i,k)|+\|f_1\|_1|A(1,k)|-|A(1,k)|-
|f_1(k)||A(1,k)|.
\end{gather*}
Using (\ref{gl1}) and (\ref{gyr}) and the inequality
$$
-1+\|f_1\|_1-|f_1(k)|\le-\mu|f_1(k)|
$$
we have
$$
\sum_{i=2}^{k-1}|A_1(i,k)|+\sum_{i=k+1}^{k+r}|A_1(i,k)|\le
\mu|A(k,k)|-\mu|f_1(k)||A(1,k))|\le\mu|A_1(k,k)|, 
$$
therefore (\ref{iroc}) holds.
\end{proof}

\begin{remark}\label{rem:schurdominant}
It is well-known that Schur complements of diagonally dominant matrices inherit the diagonal dominance property. Note that the proof of Lemma \ref{lemma:fk} also shows the corresponding result for the notion of strong diagonal dominance defined in \eqref{simplegl1}: if a matrix is diagonally dominant in the sense of \eqref{simplegl1}, for some value of $0<\mu<1$, then its Schur complements also satisfy \eqref{simplegl1} for the same value of $\mu$. 
\end{remark}

\begin{lemma}\label{lemma:gammak}
Let $A\in\mathbb{R}^{N\times N}$ be a matrix such that
\begin{equation}
\mu |A(k,k)|> \sum_{i\neq k} |A(i,k)|,\qquad k=1,\ldots,N \label{eq:hplemma}
\end{equation}
for $\mu\in (0,1)$. Since $A$ is (column) strictly diagonally dominant,
it admits an LU decomposition $A=LR$, with the usual convention that
$L$ has diagonal elements equal to 1. Then it holds
$$
|R(k,k)|\geq (1-\mu^2) |A(k,k)| ,\qquad k=1,\ldots,N.
$$ 
\end{lemma}

\begin{proof}
The thesis is obviously true for $k=1$, since $R(1,1)=A(1,1)$. Let then $k\geq 2$. It is well known from the 
theory of LU decomposition \cite{Gantmacher} that the $k$-th diagonal entry of the upper triangular factor $R$
can be written as
\begin{equation}
R(k,k)=\frac{{\rm det}\, A(1:k,1:k)}{{\rm det}\, A(1:k-1,1:k-1)}. \label{eq:detquotient}
\end{equation}
For ease of notation, denote the $k\times k$ matrix $A(1:k,1:k)$ as $\tilde{A}$.
 
The following argument is adapted from the proof of Theorem 1 in \cite{Price51}. 
Consider the problem of finding
$x_1, \ldots, x_{k-1}\in\mathbb{R}$ such that
\begin{equation}
[x_1,\, \dots\, , x_{k-1} ,\, 1]\, \tilde{A}(1:k, 1:k-1) = [0 ,\, \dots \, , 0].\label{eq:auxsystem}
\end{equation}
This amounts to solving a linear system with strictly diagonally dominant matrix $\tilde{A}(1:k-1, 1:k-1)^T$, therefore the solution exists and is unique. Moreover, it holds $|x_j|\leq \mu$, $j=1,\ldots, k-1$. Indeed, let $\ell\in \{1,\dots,k-1\}$ be such that $|x_{\ell}|=\max \{|x_1|,\ldots,|x_{k-1}|\}$ and
suppose by contradiction that $|x_{\ell}|> \mu$. The $\ell$-th equation in \eqref{eq:auxsystem} takes the form   
$$
\sum_{i=1}^{k-1} x_i \tilde{A}(i,\ell) + \tilde{A}(k,\ell)=0,
$$
which we can rewrite as
$$
-\tilde{A}(\ell,\ell) x_{\ell} = \tilde{A}(k,\ell) + \sum_{\substack{i=1 \\ i\neq\ell}}^{k-1} x_i \tilde{A}(i,\ell).
$$
Taking absolute values we have 
$$
|\tilde{A}(\ell,\ell)|\, |x_{\ell}| \leq |\tilde{A}(k,\ell)| + \sum_{\substack{i=1 \\ i\neq\ell}}^{k-1}|x_i|\, |\tilde{A}(i,\ell)| 
$$
and dividing by $|x_{\ell}|$ we obtain
$$
|\tilde{A}(\ell,\ell)| \leq \frac{|\tilde{A}(k,\ell)|}{|x_{\ell}|} + \sum_{\substack{i=1 \\ i\neq\ell}}^{k-1} |\tilde{A}(i,\ell)| \frac{|x_i|}{|x_{\ell}|}
< \frac{1}{\mu} |\tilde{A}(k,\ell)| + \sum_{\substack{i=1 \\ i\neq\ell}}^{k-1} |\tilde{A}(i,\ell)|
$$
which contradicts \eqref{eq:hplemma}.

Now, let $\hat{A}$ the matrix obtained from $\tilde{A}$ by replacing the last row with $\tilde{A}(k,:)+\sum_{i=1}^{k-1}x_i\tilde{A}(i,:)$. Note that all the entries of the last row of $\hat{A}$ are zero, except for the last one.
Therefore we have
\begin{eqnarray*}
&&{\rm det}\, \tilde{A} = {\rm det}\, \hat{A} = \hat{A}(k,k)\, {\rm det}\, \tilde{A}(1:k-1,1:k-1)=\\
&&=\left(\tilde{A}(k,k)+\sum_{i=1}^{k-1}x_i\tilde{A}(i,k)\right) {\rm det}\, \tilde{A}(1:k-1,1:k-1)
\end{eqnarray*}
and therefore, using $|x_i|\leq \mu$ and \eqref{eq:hplemma}:
\begin{eqnarray*}
|{\rm det}\, \tilde{A}|&\geq& \left(|\tilde{A}(k,k)|-\Big|\sum_{i=1}^{k-1}x_i\tilde{A}(i,k)\Big|\right) {\rm det}\, \tilde{A}(1:k-1,1:k-1)\\
&\geq& \left(|\tilde{A}(k,k)|-\mu\sum_{i=1}^{k-1}|\tilde{A}(i,k)|\right) {\rm det}\, \tilde{A}(1:k-1,1:k-1)\\
&\geq& (1-\mu^2) \,|\tilde{A}(k,k)|\, {\rm det}\, \tilde{A}(1:k-1,1:k-1).\label{eq:detmu}
\end{eqnarray*}
By substituting this inequality in \eqref{eq:detquotient} we obtain the thesis.

\end{proof}

We will also need the following result on Green generators for inverses of partially factorized matrices.

\begin{proposition}\label{prop:pfact}
Let $A=\{A(i,j)\}_{i,j=1}^N$ be a strongly regular lower band matrix of order $r$ and let $ A = L^{(\ell)} A^{(\ell)} $ be its partial LU decomposition obtained after $\ell$ steps of Gaussian elimination, with $1\leq \ell \leq N-r$. Denote as $p(i),q(i),a(i)\;(i=1,\dots,N-r)$,
$P_{N-r+1}$ the Green generators of $A^{-1}$ computed as in Theorem \ref{UF1}.

Let $\tilde{A}^{(\ell)}=A^{(\ell)}(\ell+1:N,\ell+1:N)$ and $B^{(\ell)}=(\tilde{A}^{(\ell)})^{-1}$. Then the Green generators for  $B^{(\ell)}$ computed by applying Theorem \ref{UF1} to $\tilde{A}^{(\ell)}$ are $p(i+\ell),q(i+\ell),a(i+\ell)\;(i=1,\dots,N-\ell-r)$, $P_{N-r+1}$.

\end{proposition}
 
\begin{proof}
Let $A=LR$ be the LU decomposition of $A$ and let the matrices $L_k$ be defined as in Theorem \ref{ICF1}. Analogously, let $\tilde{A}^{(\ell)}=L^{(\ell)}R^{(\ell)}$ be the LU decomposition of $\tilde{A}^{(\ell)}$ and let the matrices $L_k^{(\ell)}$ be defined as in Theorem \ref{ICF1}. Clearly it holds $L_k^{(\ell)}=L_{k+\ell}$, since the steps of the LU factorization of $\tilde{A}^{(\ell)}$ coincide with the last $N-r-\ell$ steps of the LU factorization of $A$. For the same reason, it holds 
\begin{equation}\label{eq:sameR}
R^{(\ell)}=R(\ell+1:N,\ell+1:N). 
\end{equation}
By Lemma \ref{AA}, the lower Green generators of $(L^{(\ell)})^{-1}$ are a subset of the lower Green generators of $L^{-1}$, that is,
\begin{equation}\label{eq:Lgen}
p_L^{(\ell)}(k)=p_L(k+\ell),\qquad 
q_L^{(\ell)}(k)=q_L(k+\ell),\qquad 
a_L^{(\ell)}(k)=a_L(k+\ell),
\end{equation}
for $k=1,\ldots,N-r-\ell$.

Now let us apply Theorem \ref{UF1} to $A$ and to $\tilde{A}^{(\ell)}$ and compare the Green generators for $A^{-1}$ and for $B^{(\ell)}$. From \eqref{masl} we immediately have $p^{(\ell)}(k)=p(k+\ell)$, $q^{(\ell)}(k)=q(k+\ell)$ and $a^{(\ell)}(k)=a(k+\ell)$, for $k=1,\ldots,N-r-\ell$. Then from \eqref{mishln} and \eqref{mishl} we obtain 
\begin{eqnarray*}
&&p^{(\ell)}(k)=\frac1{\gamma^{(\ell)}_k}(p^{(\ell)}_L(k)-X^{(\ell)}_k P^{(\ell)}_{k+1}a^{(\ell)}(k))=\\
&& = \frac1{\gamma_{k+\ell}}(p_L(k+\ell)-X_{k+\ell} P_{k+\ell+1}a(k+\ell))= p(k+\ell),
\end{eqnarray*}
where the relations $\gamma^{(\ell)}_k=\gamma_{k+\ell}$ and $X^{(\ell)}_k = X_{k+\ell}$ are a consequence of \eqref{eq:sameR}. 
\end{proof}

Inverses of strictly diagonally dominant matrices can be bounded in norm, in terms of a parameter that quantifies the amount of diagonal dominance. Indeed, recall the following result from \cite{Varah75}:
\begin{theorem}\label{thm:varah}
Let $A\in\mathbb{R}^{N\times N}$ be diagonally dominant by columns and set
$$
\beta=\min_k \left(|A(k,k)|-\sum_{j\neq k}|A(j,k)|\right).
$$
Then it holds
$$
\|A^{-1}\|_1\leq\frac{1}{\beta}.
$$
\end{theorem}

From Theorem \ref{thm:varah} we deduce the following result, which will be useful in the proof of Theorem \ref{TTR}.

\begin{cor}\label{cor:varah}
Let $A=\{A(i,j)\}_{i,j=1}^N$ satisfy the condition \eqref{gl1}. Then it holds
$$
\|A^{-1}\|_1\leq\frac{1}{(1-\mu)\min_i |A(i,i)|}.
$$
\end{cor}
\begin{proof}
Let us apply Theorem \ref{thm:varah}. In our case we have 
$$
|A(k,k)|-\sum_{j\neq k}|A(j,k)|>(1-\mu)|A(k,k)|
$$
for $k=1,\ldots,N$ and therefore
$$
\beta=\min_k \left(|A(k,k)|-\sum_{j\neq k}|A(j,k)|\right)>(1-\mu)\min_k|A(k,k)|,
$$
from which we deduce
$$
\|A^{-1}\|_1\leq\frac{1}{\beta}<\frac{1}{(1-\mu)\min_k|A(k,k)|}.
$$
\end{proof}

We are now ready to prove Theorem \ref{TTR}.

\vspace{20pt}

{\bf Proof of Theorem \ref{TTR}.}
From Theorem \ref{UF1} and Remark \ref{rem:green2scalar} we have for the Green block form of $A^{-1}$:
$$
(A^{-1})'(i,j)=p(i)a^>_{ij}q(j),\qquad 1\leq i\leq N,\quad 0\le j\le N-r
$$
with the lower Green generators $a(k),q(k)$ obtained from the lower triangular
matrix $L^{-1}$ via
$$
q(k)=q_L(k),\;a(k)=a_L(k),\;k=1,\dots,N-r
$$
and $p(i)$ obtained in Step 2 of the algorithm in Theorem \ref{UF1}.
We will prove that the generators $p(i)$ are bounded, i.e.,
\begin{equation}\label{lrlo}
\|p(i)\|_{\infty}\leq M,\quad  i=1,2,\dots N,
\end{equation}
where $M=\frac{1}{(1-\mu)\min_i |A_{i,i}|}$,
and the vectors $a^>_{ij}q(j)$ satisfy the estimates
\begin{equation}\label{lrlos}
\|a^>_{ij}q(j)\|_1\leq\gamma^{i-j-r},\quad i\geq j+r
\end{equation}
for $0<\gamma=\mu^{1/r}<1$. This will yield the thesis.

Unless otherwise noted, in this proof all vector and matrix norms are understood to be the vector 1-norm and the corresponding induced matrix norm. 

Note that from Lemma \ref{AA} we obtain the lower Green generators of the matrix 
$L^{-1}$ via the formulas (\ref{natas}) with the partitions (\ref{iof5}). Moreover, Lemma
\ref{lemma:fk} ensures $\|f_k\|_1\le\mu<1$. 

{\bf Part 1}: let us prove (\ref{lrlos}). At first we derive an auxiliary
relation. Set 
\be\label{larkr}
F_m=I_r,\quad F_i=a(i)a(i-1)\cdots a(m+1),\;i=m+1,\dots,m+r.
\end{equation}
From (\ref{natas}) we can write 
\be\label{itr}
a_L(k)=a(k)=-f_k e_1^T+J
\end{equation}
where $e_1=[1,0,\ldots,0]^T\in\mathbb{R}^r$ and $J$ is the $r\times r$ Jordan block 
$$
J=\left(\ba{cccccc}
0&1&0&\dots&0&0\\
0&0&1&\dots&0&0\\
\vdots&\vdots&\vdots&\vdots&\vdots&\vdots\\
0&0&0&\dots&0&1\\0&0&0&\dots&0&0\ea\right).
$$
Note that \eqref{itr} implies $\|a(k)\|\leq 1$.

Set 
\be\label{olg}
V_i=J^{i-m},\;i=m,\dots,m+r.
\end{equation}
It is clear that $V_m=I_r,\;V_{i+1}=JV_i,\;i=m,\dots,m+r-1,\;V_{m+r}=0$.
Next we set $\Gamma_i=F_i-V_i,\;i=m,\dots,m+r$, i.e.
\be\label{lariss}
F_i=\Gamma_i+V_i,\quad i=m,m+1,\dots,m+r.
\end{equation}
Moreover we have $F_{i+1}=a(i+1)F_i$ and using (\ref{itr}) and (\ref{lariss}) we
get
\be\label{irg}
F_{i+1}=(-f_{i+1}e_1^T+J)(\Gamma_i+V_i)=
-f_{i+1}e_1^T\Gamma_i+J\Gamma_i-f_{i+1}e_1^T\Gamma_i+V_{i+1},\;i=m,\dots,m+r-1.
\end{equation}
We have obviously $\Gamma_m=F_m-V_m=0$ and moreover (\ref{irg}) and 
(\ref{lariss}) imply
\be\label{irut}
\Gamma_{i+1}=-f_{i+1}e_1^T\Gamma_i+J\Gamma_i-f_{i+1}e_1^TV_i,\;i=m,\dots,m+r.
\end{equation}
Hence for any $r$-dimensional column vector $g=[g(1),\dots,g(r)]^T$ we get
\be\label{alk}
\Gamma_{i+1}g=-f_{i+1}e_1^T(\Gamma_ig)+J\Gamma_ig-f_{i+1}e_1^TV_ig.
\end{equation}
We prove by induction that
\begin{equation}
\|\Gamma_ig\|\le\mu(|g(1)|+\dots+|g(i-m)|),\;i=m,\dots,m+r.\label{alu1}
\end{equation}
For $i=m$ since $\Gamma_m=0$ the statement is clear. Assume that for some $i$ 
with $m\le i\le m+r-1$ the relation (\ref{alu}) holds. Using (\ref{alk}) we get
$$
\|\Gamma_{i+1}g\|\le\|f_{i+1}\|\, |e^T_1(\Gamma_ig)+\|J(\Gamma_ig)\|+
\|f_{i+1}\|\, |g(i-m+1)|.
$$
Using the inequalities $\|f_{i+1}||\le \mu<1$ we get
$$
\|\Gamma_{i+1}g\|\le |e^T_1(\Gamma_ig)|+\|J(\Gamma_ig)\|+\mu|g(i-m+1)|.
$$
For any $r$-dimensional column vector $h$ the equality 
$$
|e_1^Th|+\|Jh\|=\|h\|
$$
holds and using (\ref{alu1}) we conclude that 
$$
\|\Gamma_{i+1}g\|\le\mu(|g(1)|+\dots+|g(i-m)|+|g(i-m+1)|).
$$
Taking $i=m+r$ in (\ref{alu1}) we get $\|\Gamma_{m+r}\|\le\mu$.

Moreover the equality (\ref{olg}) 
with $i=m+r$ implies $V_{m+r}=0$ and therefore $F_{m+r}=\Gamma_{m+r}$ and we 
conclude that
\be\label{lij}
\|F_{m+r}\|\le\mu.
\end{equation}

In accordance with (\ref{lrlos}) consider 
the products $a_{ij}^>,\;i-j\geq r$. Recall that $a_{ij}^>$ is defined as the product of $i-j-1$ factors.  We use the representation $i-j-1=rt+l$, where the 
positive integer $t$ is the quotient of $i-j-1$ divided by $r$, and $l$ with $0\le l<r$ is the reminder. We have
$$
a^>_{ij}=a(i-1)\cdots a(j+1)=a(i-1)\cdots a(i-l)a(rt+j)\cdots a(j+1),
$$
where  
the product $a(i-1)\cdots a(i-l)$ is assumed to be $1$ if $l=0$.
From the inequalities $\|a(k)\|\le1$ applied with $k=i-l,\dots,i-1$ and from (\ref{lij}) we get
$$
\|a_{ij}^>\|\le \|a(rt+j)\cdots a(j+1)\|\le\mu^t.
$$
Using the equality $t=(i-j)/r-(l+1)/r$ and setting 
$\gamma=\mu^{1/r}$ we conclude that
\be\label{alu}
\|a_{ij}^>\|\le \gamma^{i-j-r}.
\end{equation}
From 
\eqref{natas}, \eqref{iof5} and \eqref{masl} 
we deduce $\|q(k)\|=1$,
whichi concludes the proof of \eqref{lrlos}.

{\bf Part 2}: 
Consider at first the case $i=1$. From the Green representation of $A^{-1}$ we have 
$p(1,:)=A^{-1}(1,1:r)$. But Corollary \ref{cor:varah} implies 
$$|A^{-1}(j,k)|\leq \frac{1}{(1-\mu)\min_i |A_{i,i}|}<M$$ for all indices $j,k=1,\ldots,N$. Therefore it holds
$$\|p(1,:)\|_{\infty}\leq M.$$

Now suppose $i>1$ and consider the partial LU decomposition 
$$ A = L^{(i)} A^{(i)} $$
obtained after $i-1$ steps of Gaussian elimination as in Theorem \ref{ICF1}. Note that the matrix $A^{(i)}(i:N,i:N)\in\mathbb{R}^{(N-i+1)\times (N-i+1)}$ inherits from $A$ the diagonal dominance property \eqref{gl1}, as proved in Lemma \ref{lemma:fk}; see also Remark \ref{rem:schurdominant}.

Let $\tilde{A}^{(i)}=A^{(i)}(i:N,i:N)$ and $B^{(i)}=(\tilde{A}^{(i)})^{-1}$. From Proposition \ref{prop:pfact} we have $p(i,:)=B^{(i)}(1,1:r)$. The definition of matrix 1-norm implies that the $1$-norm of any given matrix is larger than or equal to the absolute value of any entry in the matrix; in particular it holds $|B^{(i)}(h,k)|\leq \|B^{(i)}\|_1$ for any index pair $(h,k)$ and therefore $\|B^{(i)}(1,1:r)\|_{\infty}\leq \|B^{(i)}\|_1$. Then Corollary \ref{cor:varah} applied to $\tilde{A}^{(i)}$ yields  
\begin{equation}\label{eq:firstestimate}
\|p(i,:)\|_{\infty}\leq\frac{1}{(1-\mu)\min_j|\tilde{A}^{(i)}(j,j)|}.
\end{equation}
In order to give an estimate for $\tilde{A}^{(i)}(j,j)$, observe that the upper triangular factor in the LU factorization of $\tilde{A}^{(i)}(j,j)$ coincides with $R(i:N,i:N)$, where $R$ is the upper triangular factor in the LU factorization of $A$. From Lemma \ref{lemma:gammak} and from the diagonal dominance property of $\tilde{A}^{(i)}$ we have
$$
(1+\mu^2)|\tilde{A}^{(i)}(j,j)|\geq |R(i+j-1,i+j-1)|\geq (1-\mu^2)|A(i+j-1,i+j-1)|,
$$
and therefore
\begin{equation}\label{eq:diagbound}
\tilde{A}^{(i)}(j,j)|\geq \frac{1-\mu^2}{1+\mu^2}|A(i+j-1,i+j-1)|.
\end{equation}
From \eqref{eq:firstestimate} and \eqref{eq:diagbound} we deduce 
$$\|p(i,:)\|_{\infty}\leq\frac{1+\mu^2}{(1-\mu)(1-\mu^2)\min_k|A(k,k)|},$$
which concludes the proof.

$\hfill\Box$

%\textcolor{red}{Check norms, 1 or max.}

\section{QR-based decay bounds}\label{sec:QR}
The decay bounds presented so far rely on the LU-based inversion algorithm developed in \cite{BE23}. The same paper includes an inversion algorithm for banded matrices that relies on the QR factorization. This can also be used to deduce decay bounds for the inverse of a banded matrix, in combination with a suitable diagonal dominance condition. 
With this approach, however, the proof is more involved with respect to the LU-based formulas and the resulting bounds are less effective. We report here these bounds for completeness; the proof is omitted. 

\begin{theorem}\label{TTRQR}
Let $A=\{A_{ij}\}_{i,j=1}^N$ be an invertible lower-banded matrix of order $r$. 
Assume that 
$$
\sum_{i=1}^{k-1}\sum_{j=k+1}^N|A(i,j))|^2\le C_0,\quad k=1,\dots,N-r
$$
for a suitable constant $C_0>0$, which does not depend on $N$ or $k$. 

Assume also that $A$ satisfies the following strong dominance condition:
$$
|A(k,k)|\ge K\sqrt{\sum_{i=1}^{k-1}|A(i,k)|^2+\sum_{i=k+1}^{k+r}|A(i,k)|^2}+1
$$
with $K$ such that
$$
K\ge \max\left\{\frac2{x_0},4(3+2C_0r\sqrt{r}),
2\sqrt{r^3\left(\frac{\sqrt{3}+1}{2}\right)^{2r}-1}\right\}. 
$$
Then there exist $M>0$ and $\gamma$ with $0\le\gamma<1$ such that
$$
|(A^{-1})(i,j)|\le M\gamma^{i-j},\quad i\ge j.
$$
Moreover, $M$ and $\gamma$ can be explicitly computed as
$$
M = 2\mu+1,\qquad \gamma=(\mu r\sqrt{r})^{1/r},
$$
where
$$
\mu=\frac{\delta}{\sqrt{1+\delta^2}},\qquad \delta = \frac{2}{K}.
$$
\end{theorem}

\begin{remark} In the case of a two-sided banded matrix $A$ of order $r$, the diagonal dominance 
conditions in Theorem \ref{TTRQR} take the form
$$
\sum_{i=k-r}^{k-1}\sum_{j=k+1}^{k+r}|A(i,j))|^2\le C_0
$$
and
$$
|A(k,k)|\ge K\sqrt{\sum_{i=k-r}^{k-1}|A(i,k)|^2+\sum_{i=k+1}^{k+r}|A(i,k)|^2}+1.
$$
\end{remark}

QR-based bounds from Theorem \ref{TTRQR} share some of the qualitative features of LU-based ones from Theorem \ref{TTR}: they  are easily computable and can be applied to classes of nonsymmetric matrices for which usual approximation-based bounds are not effective. However, they are typically more pessimistic than LU-based bounds, as shown in Example 4. %Moreover, the constraints on $K$ imposed by Theorem \ref{TTRQR} can be rather restrictive, as shown...

\section{Numerical experiments}\label{sec:numerical}
% All experiments to be re-done.

In this section we investigate experimentally the behavior of the bounds from Theorem \ref{thm:bound} and we compare them to other bounds available in the literature, with the goal of illustrating their strong and weak points. In particular, we compare our bounds, when possible, to results from \cite{DMS84} (denoted in figure legends as DMS), from \cite{Frommer18, FrommerNLAA} and from \cite{CH83, Hasson07}.

Unless otherwise noted, in all the experiments we have chosen $\mu$ as the smallest value that satisfies \eqref{gl1}. By ``exact inverse'' of a matrix we mean the inverse computed by the MATLAB command {\tt inv}.

The MATLAB code used to generate the examples in this section is available from {\tt https://people.cs.dm.unipi.it/boito/QSdecay.zip}.

\vspace{15pt}

{\bf Example 1: real symmetric, two-sided banded matrices.} 

{\bf (a)} Let us start with a test on a symmetric positive definite, Toeplitz, two-sided banded matrix of size $50\times 50$ with bandwidth $r=3$:
\begin{equation}
A(i,j)=\left\{\begin{array}{lll}
6.25 & {\rm if} & i=j,\\
0.25 & {\rm if} & i\neq j,\, |i-j|\leq 3,\\
0 & {\rm otherwise.}\\
\end{array}\right.\label{eq:Aex1}
\end{equation}
We compare the bounds from Theorem \ref{thm:bound} with the bounds from \cite{DMS84}. Recall that the work by Demko, Moss and Smith \cite{DMS84} relies on polynomial approximation of the function $f(x)=1/x$ on a real interval $[a,b]$ with $0<a<b$ and, as such, requires spectral information on $A$. Ideally one should choose $a$ and $b$ as the minimum and maximum eigenvalue of $A$, respectively; if these are not available, one may choose instead lower/upper bounds on the spectrum, such as the estimates given by Gershgorin's theorems, although the resulting bounds will be less sharp. 

\begin{figure}
\begin{center}
\includegraphics[width=0.45\textwidth]{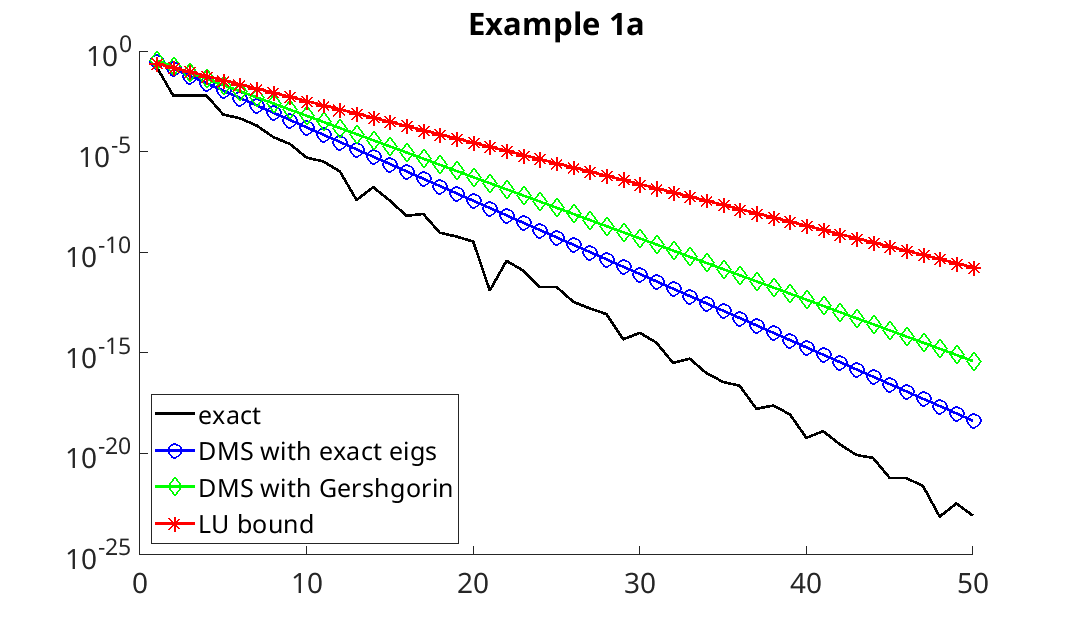}
\includegraphics[width=0.45\textwidth]{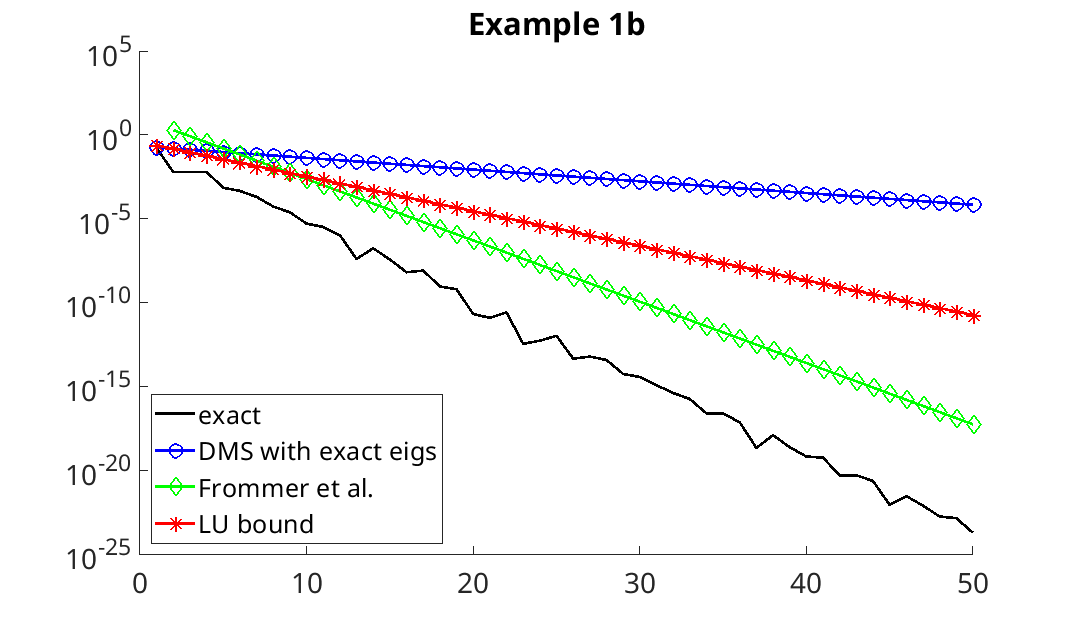}

\includegraphics[width=0.45\textwidth]{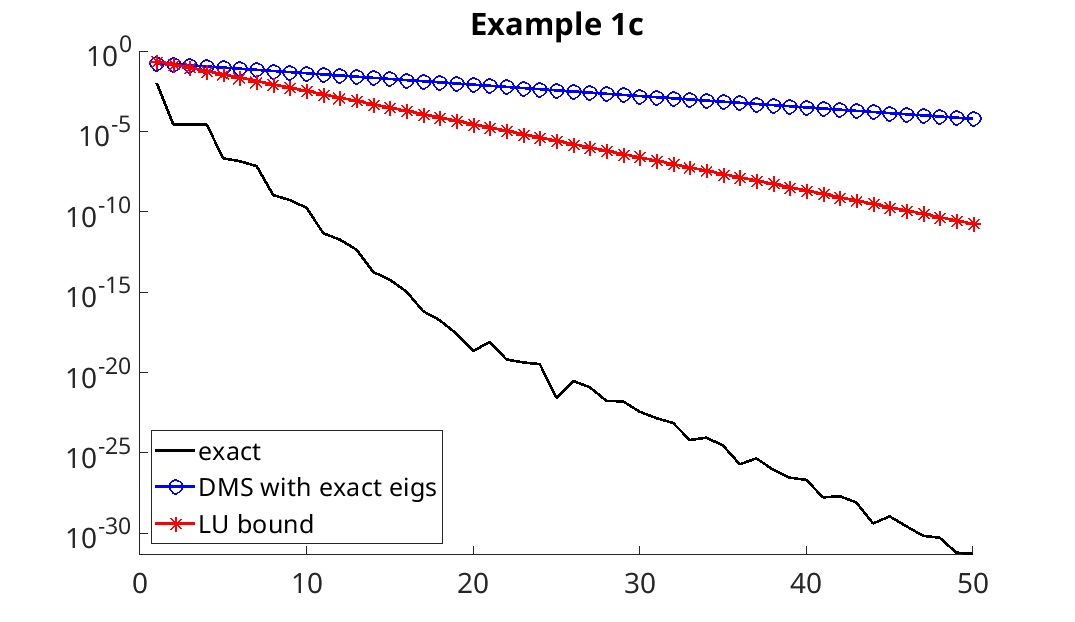}
\includegraphics[width=0.45\textwidth]{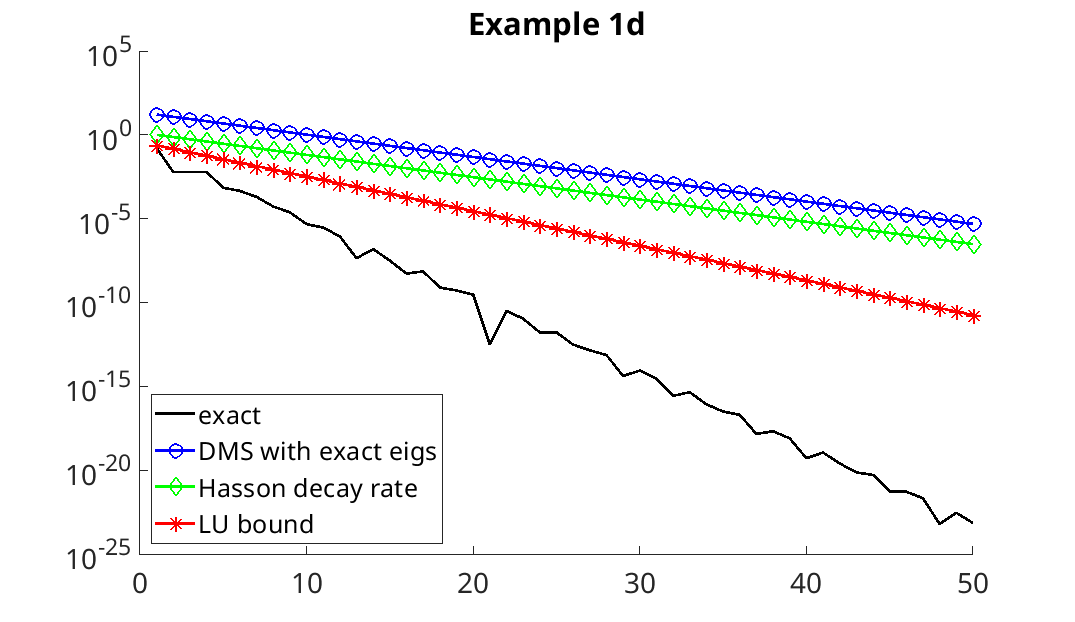}
\caption{\footnotesize Plots for Example 1: behavior of the absolute value of the first column of $A^{-1}$, in logarithmic scale. Row index is on the $x$ axis, entry values and bounds on the $y$ axis. DMS stands for the bounds from \cite{DMS84}, LU bounds are from Theorem \ref{thm:bound}.}\label{fig:ex1}
\end{center}
\end{figure}

The upper left plot in Figure \ref{fig:ex1} clearly shows that our bounds are more pessimistic than the ones in \cite{DMS84}. This is to be expected, because the decay rate from \cite{DMS84} is known to be optimal for Toeplitz, symmetric positive definite matrices. More generally, we do not expect our bounds to be competitive for symmetric, positive definite matrices.

{\bf (b)} Let us now introduce a single large eigenvalue by setting $A(20,20)=100$. Bounds from \cite{DMS84} are now quite pessimistic. In fact, this is a good test matrix for bounds from \cite{Frommer18}, Theorem 3.1, which provide a very good decay estimate relying on the notion of ``effective condition number''. On this example, these bounds yield
\begin{equation}
|[A^{-1}]_{i,j}|\leq C q_1^{\frac{|i-j|}{r}-1}, \qquad C=\frac{2}{\lambda_1},
\quad q_1=\frac{\sqrt{\kappa_e(A)}-1}{\sqrt{\kappa_e(A)}+1},\quad \kappa_e(A)=\frac{\lambda_{N-1}}{\lambda_1},\label{eq:frommer}
\end{equation}
where $\lambda_1\leq\ldots\lambda_{N-1}\leq\lambda_N$ are the eigenvalues of $A$, and $|i-j|\geq r$. The performance of our LU-based decay bounds is between these two, as illustrated by the upper right plot in Figure \ref{fig:ex1}.

{\bf (c)} Suppose now that half of the eigenvalues are large and clustered. To this end, we take the matrix from \eqref{eq:Aex1} and add $100$ to the first $25$ diagonal entries. The bounds \eqref{eq:frommer} cannot be applied in this case. LU-based bounds perform better than the bounds from \cite{DMS84}, although they are still somewhat pessimistic: see the lower left plot in Figure \ref{fig:ex1}.

{\bf (d)} Let us modify our matrix so that it is symmetric indefinite. We take again $A$ from \eqref{eq:Aex1} and change the sign of a few diagonal entries, e.g., for indices 10, 11 and 12. The paper \cite{DMS84} provides bounds also for this case, relying essentially on the application of polynomial approximation of $f(x)=1/x$ to the spectrum of $A^TA$. Figure \ref{fig:ex1} (lower right plot) shows that now the bounds from Theorem \ref{thm:bound} perform better than \cite{DMS84}. Bounds from Corollary 2 in \cite{FrommerNLAA} give in this case the same decay rate as in \cite{DMS84}.  

For completeness we report also the decay rate implied by Chui and Hasson's results on polynomial approximation of $f(x)=1/x$ on disjoint real intervals \cite{CH83, Hasson07}. Given $a,b\in\mathbb{R}$ with $0<a<b$, there exists a constant $C$ such that for any symmetric, $r$-banded matrix $A$ with eigenvalues contained in $[-b,-a]\cup [a,b]$, it holds
\begin{equation}\label{eq:hasson}
|[A^{-1}]_{ij}|\leq C \lambda ^n,\qquad \lambda=\left(\frac{b-a}{b+a}\right)^{1/2r}.
\end{equation}
However, these results do not yield an actual decay bound, since $C$ is not known explicitly. The decay rate is again less favorable than the one from Theorem \ref{thm:bound}.  

\vspace{20pt}

In fact, a comparison between our bounds and the bounds from \cite{DMS84} with Gershgorin estimates can also be done analytically for specific classes of matrices, including the ones used in items (a) and (d) in this example. If $A$ is a two-sided $r$-banded real symmetric positive definite matrix, then the decay rate for its inverse (i.e., the constant $\xi$ in \eqref{generalbound}) proposed in \cite{DMS84} is
$$
\lambda_0 = \left(\frac{\sqrt{b/a}-1}{\sqrt{b/a}+1}\right)^{1/r},
$$
where $0<a<b$ are such that the spectrum of $A$ is contained in $[a,b]$. Assume that $A$ satisfies \eqref{gl1} and has constant positive diagonal. Gershgorin's theorems imply that we can take $a = 1-\mu$ and $b=1+\mu$. With this choice we have
$$
\lambda_0=\left(\frac{1-\sqrt{1-\mu^2}}{\mu}\right)^{1/r},
$$
which is easily seen to be smaller than the decay rate $\gamma=\mu^{1/r}$ from Theorem \ref{thm:bound}. 

On the other hand, if $A$ is two-sided $r$-banded real symmetric indefinite, then the decay rate from \cite{DMS84} is
$$
\lambda_1 = \left(\frac{b/a-1}{b/a+1}\right)^{1/2r},
$$
where $0<a<b$ are such that the spectrum of $A$ is contained in $[-b,-a]\cup [a,b]$. Assume that all diagonal entries of $A$ are equal in absolute value. Then we can take again $a = 1-\mu$ and $b=1+\mu$ and we obtain
$\lambda_1=(\sqrt{\mu})^{1/r}$,
that is, a worse decay rate than the one from Theorem \ref{thm:bound}. A similar conclusion holds for the decay rate \eqref{eq:hasson}, which in this case turns out to be equal to $\lambda_1$.

\vspace{20pt}

{\bf Example 2: a nonsymmetric, two-sided banded matrix with a large eigenvalue.} 

Let us take $A$ as in Example 1a and modify it as follows:

{\tt A(21:23,20) = A(21:23,20)*25; A(20,20) = -100;} 

{\tt A(22:24,21) = A(22:24,21)*25; A(21,21) = -100;} 
 
{\tt A(1:33,30) = A(1:33,30)/100; A(30,30) = 1;} 

Note that $A$ has two complex conjugate eigenvalues close to $-100$, one eigenvalue close to $1$ and all the other eigenvalues in the interval $[5.6,7.7]$. As expected, this is a case where the bounds from \cite{DMS84} do not perform well.

\begin{figure}
\begin{center}
\includegraphics[width=0.8\textwidth]{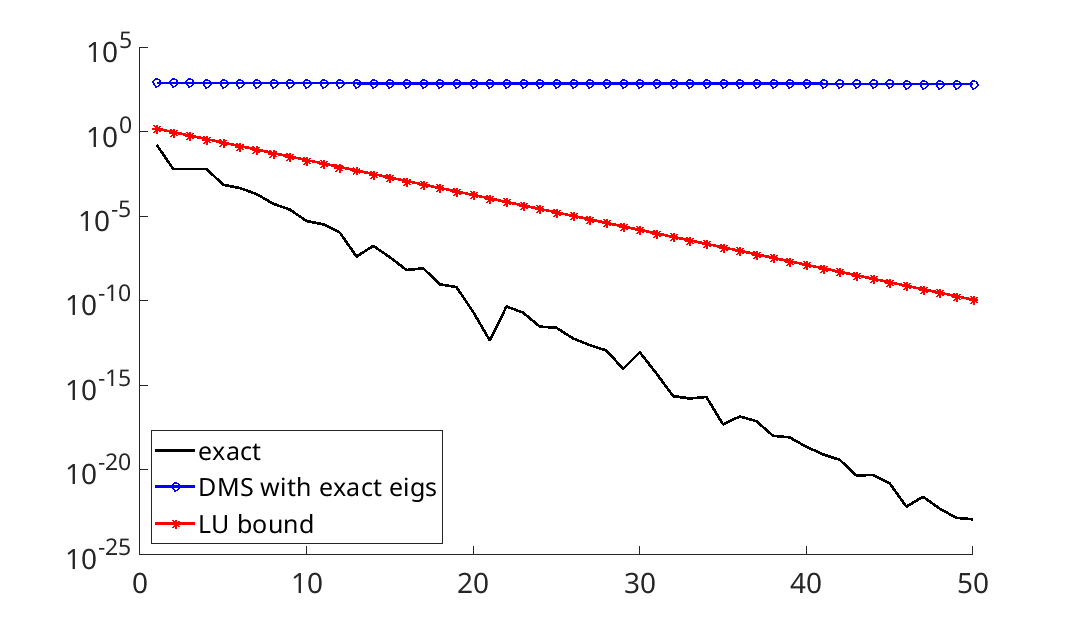}
\caption{\footnotesize LU bounds from Theorem \ref{thm:bound} for Example 2. As usual, bounds are compared to the behavior of the first column of $A^{-1}$. Row index is on the horizontal axis. DMS stands for the bounds from \cite{DMS84}.}\label{fig:ex2}
\end{center}
\end{figure}

% Normal matrices? Thm 2 from Frommer et al, NLAA 2018, also try test on their shifted skew-Hermitian matrix. 

\vspace{20pt}

\begin{figure}
\begin{center}
\includegraphics[width=0.45\textwidth]{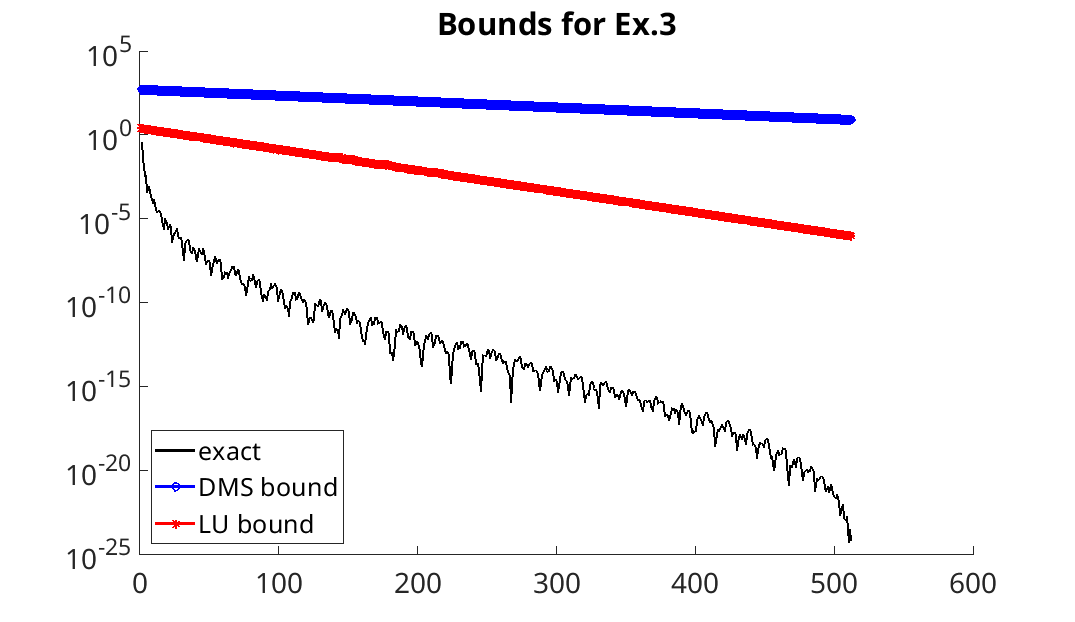}
\includegraphics[width=0.45\textwidth]{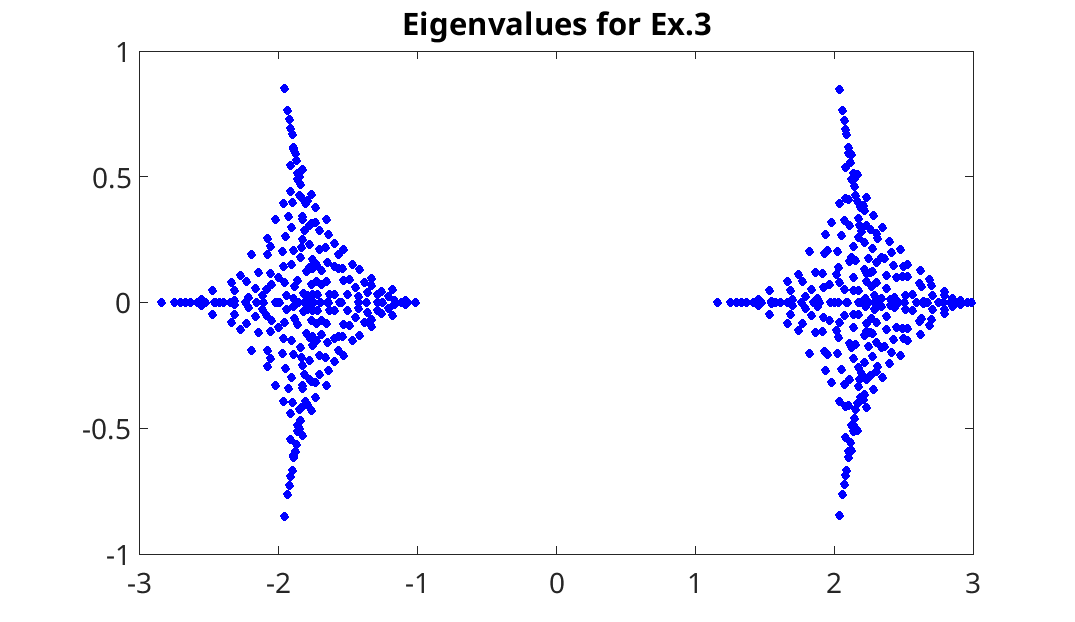}
\caption{\footnotesize Left plot: decay bounds for Example 3. As usual, bounds are compared to the behavior of the first column of $A^{-1}$; row index is on the horizontal axis. DMS stands for the bounds from \cite{DMS84}. Right plot: eigenvalues of matrix $A$.}\label{fig:gre}
\end{center}
\end{figure}

{\bf Example 3: nonsymmetric, two-sided banded matrices with clustered eigenvalues.}
Here $A$ is the $512\times 512$ matrix {\tt gre\_512} from the Sparse Matrix Collection \cite{SMC}, plus a diagonal term $I_{256}\oplus (-I_{256})$. The eigenvalues form two clusters, one with positive real part and one with negative real part; see right plot in Figure \ref{fig:gre}. The matrix has lower bandwidth $r=24$ and upper bandwidth $64$. Figure \ref{fig:gre} (left plot) compares bounds from Theorem \ref{thm:bound} and from \cite{DMS84}.

\vspace{20pt}

{\bf Example 4: one-sided banded matrices.} 
In this example we choose $A$ as a $100\times 100$ lower-banded matrix with bandwidth $5$ and eigenvalues both in the left and right complex half plane. Although not optimal, our bounds manage to capture the exponential decay in the lower triangular part of $A^{-1}$. We only compare our bounds to the exact inverse, as we are not aware of other explicit bounds in the literature that apply to such examples. 

\begin{figure}
\begin{center}
\includegraphics[width=0.45\textwidth]{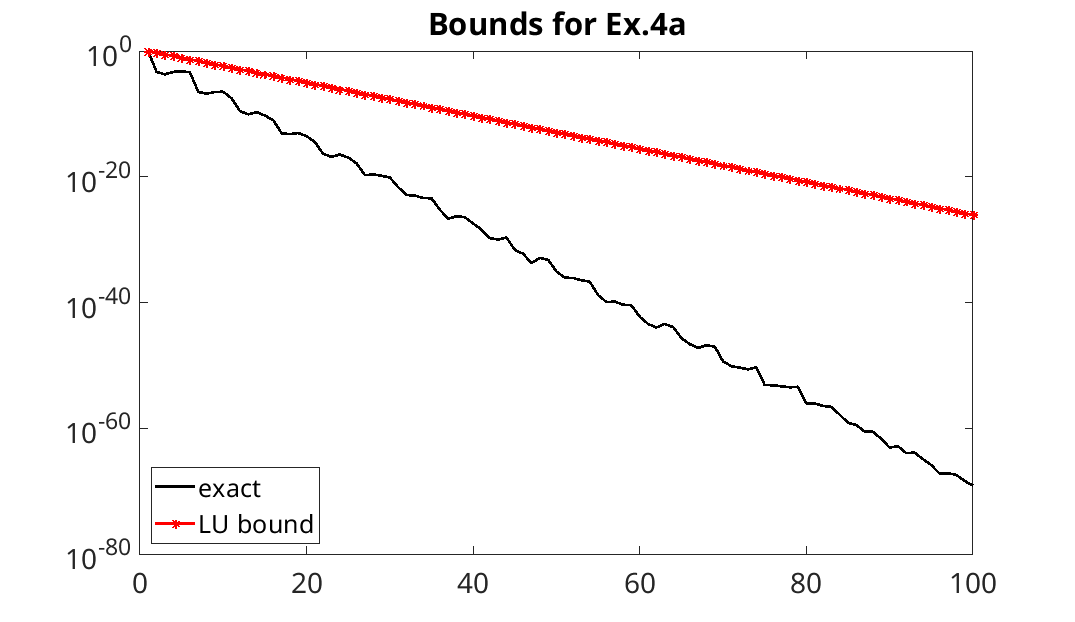}
\includegraphics[width=0.45\textwidth]{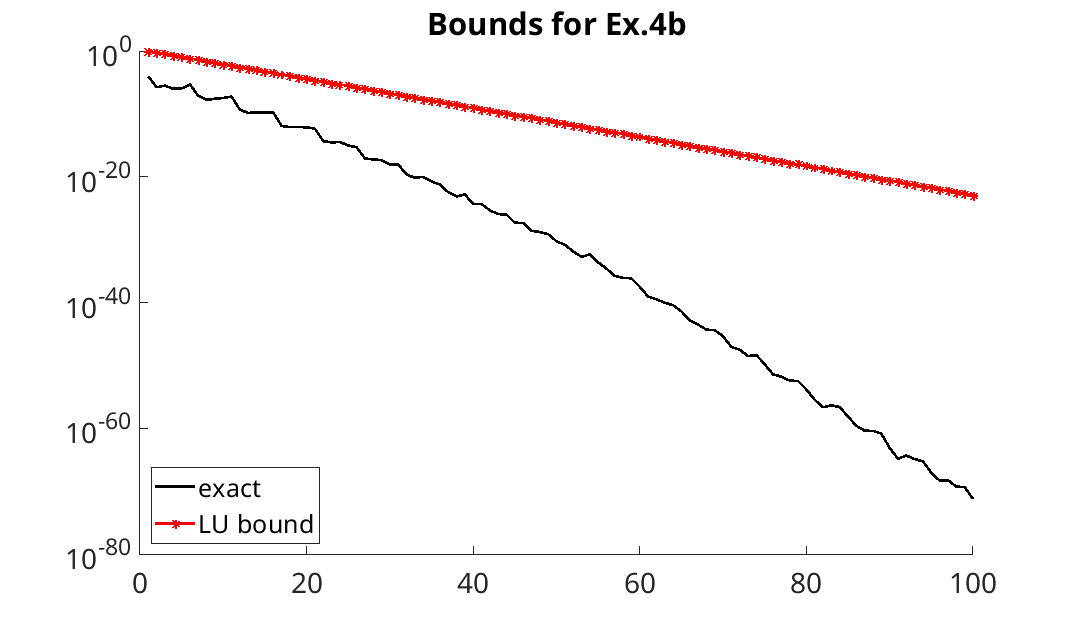}

\includegraphics[width=0.45\textwidth]{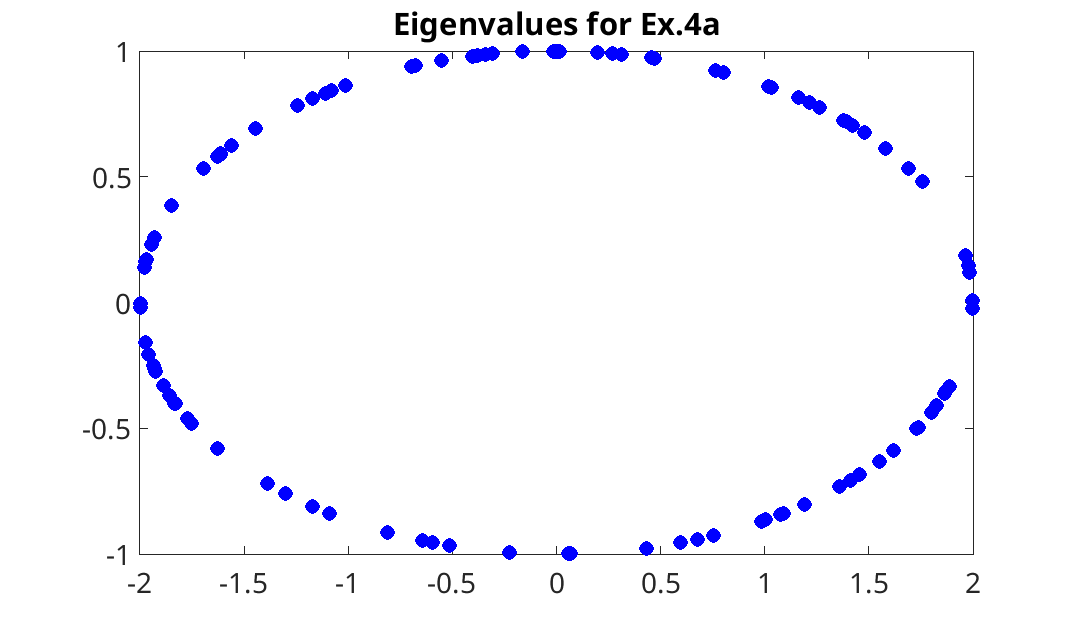}
\includegraphics[width=0.45\textwidth]{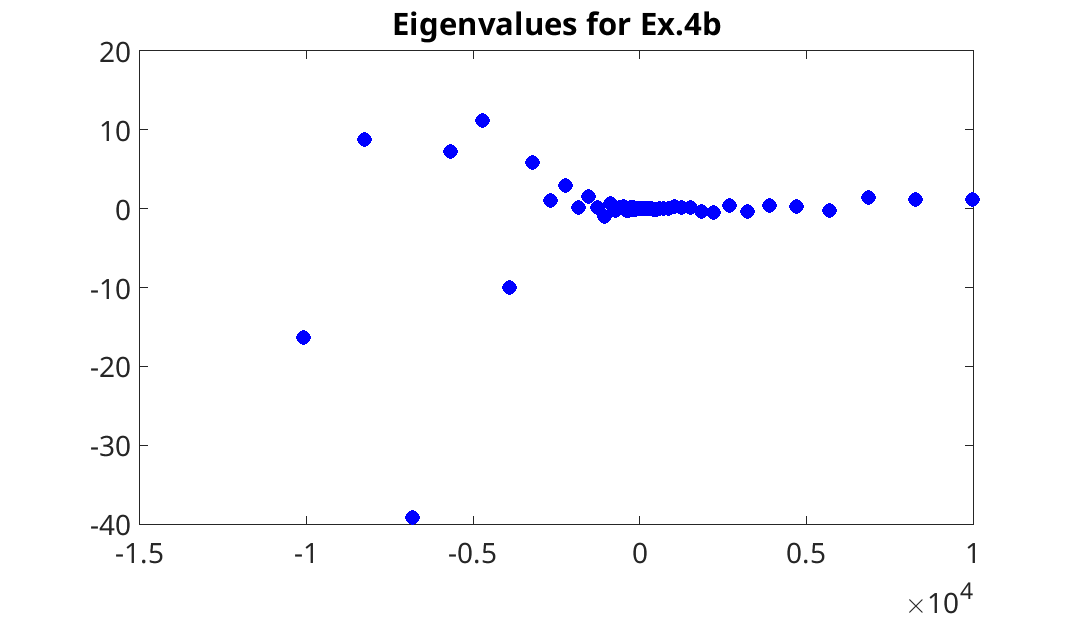}
\caption{\footnotesize First plot: LU bounds for Example 4a. Second plot: LU bounds for Example 4b. As usual, bounds are compared to the behavior of the first column of $A^{-1}$. Row index is on the horizontal axis. The two bottom plots show eigenvalue distributions in the complex plane for the two examples.}\label{fig:onesided}
\end{center}
\end{figure}

{\bf (a)} The eigenvalues of $A$ are randomly chosen close to an ellipse centered in $0$, with semiaxis 2 along the real axis and 1 along the imaginary axis. Nonzero, non-diagonal entries are random in $[-10^{-3},10^{-3}]$, with uniform distribution. See Figure \ref{fig:onesided}, left plots. 

{\bf (b)} The eigenvalues of $A$ are such that their real parts roughly follow a logarithmic distribution in $[-10^4,-1]\cup [1,10^4]$. See Figure \ref{fig:onesided}, right plots. 

\vspace{20pt}

\begin{figure}
\begin{center}
\includegraphics[width=0.45\textwidth]{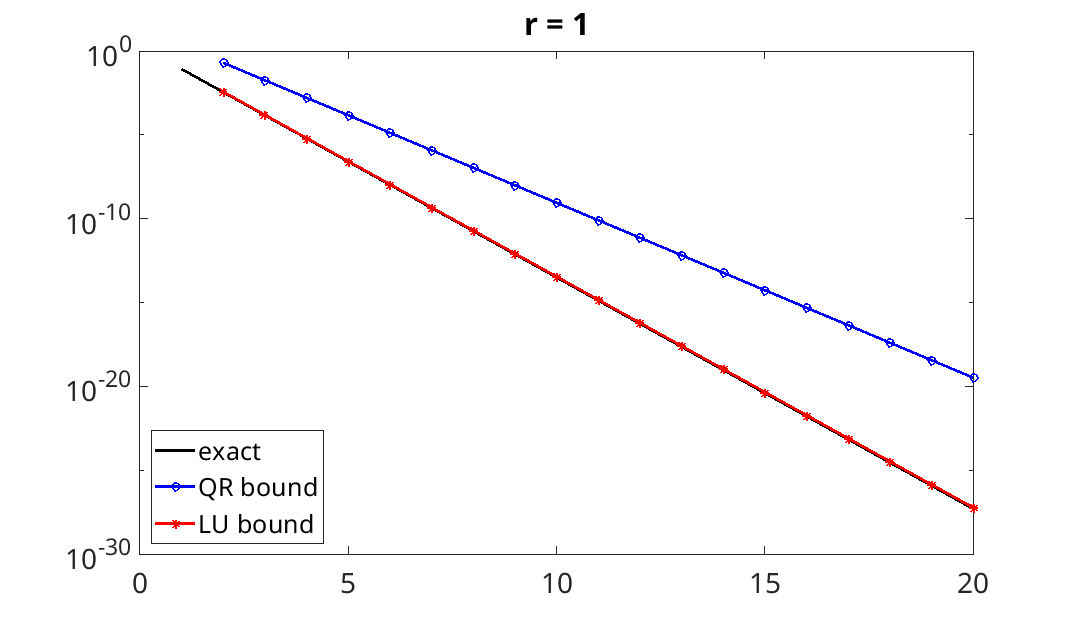}
\includegraphics[width=0.45\textwidth]{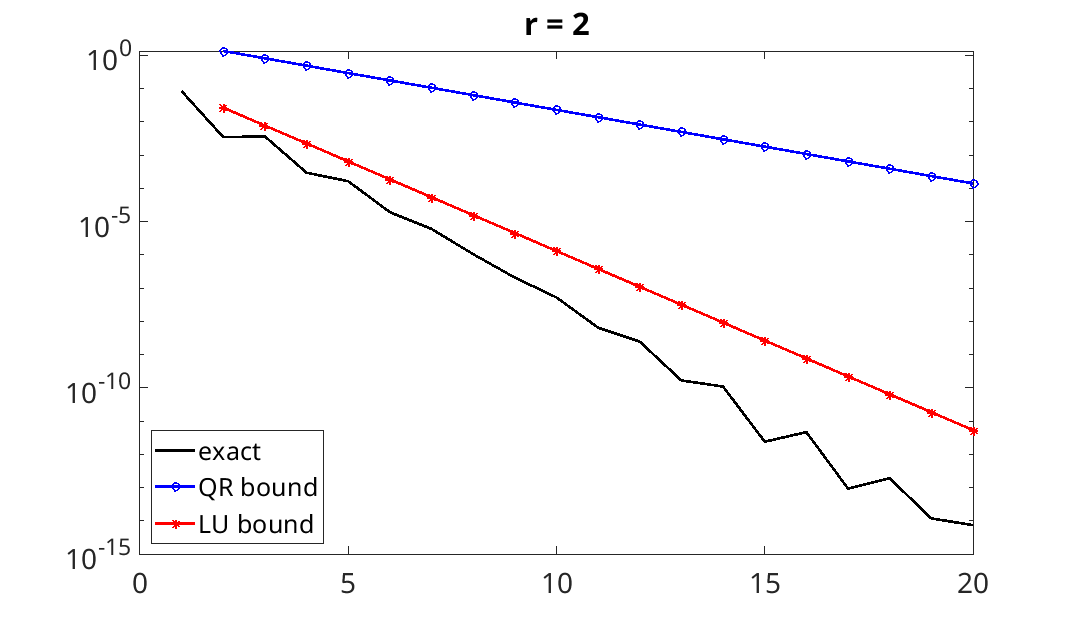}
\caption{\footnotesize Comparison between LU- and QR-based bounds from Theorems \ref{thm:bound} and \ref{TTRQR} (Example 5). As usual, bounds are compared to the behavior of the first column of $A^{-1}$. Row index is on the horizontal axis. In the left plot, LU-based bounds are very close to exact values.}\label{fig:LUvsQR}
\end{center}
\end{figure}

{\bf Example 5: comparison between LU- and QR-based bounds.} 
Here $A$ is a $20\times 20$ lower-banded matrix with nonzero subdiagonal entries equal to $0.5$ and exponentially decaying upper triangular part. Eigenvalues are complex and clustered around $\pm 12$. Figure \ref{fig:LUvsQR} compares bounds from Theorems \ref{thm:bound} and \ref{TTRQR} for two cases: lower bandwidth $r=1$, i.e., Hessenberg structure (left plot), and  $r=2$ (right plot). Clearly, LU-based bounds are more accurate, particularly for larger bandwidth.   

\vspace{20pt}

\section{Conclusions}\label{sec:conclusions}

We have presented new exponential decay bounds for inverses of banded matrices under hypotheses of strong diagonal dominance. In our opinion, these results have several appealing features:
\begin{itemize}
\item ease of computation, as these bounds do not require spectral information or computation of the field of values;
\item they can be applied to classes of matrices that do not satisfy the hypotheses required by approximation-based bounds; 
\item they rely on a quasiseparable approach that has not been exploited so far in the literature, except for particular cases. Even from a purely theoretical perspective, we think it is interesting {\em per se} that exponential decay bounds can be obtained also in this way.
\end{itemize}  

%The bounds are likely not optimal, as illustrated by tests on symmetric positive definite matrices, but they still outperform previously known bounds for some classes of indefinite or nonsymmetric matrices. 

A natural question is whether the approach presented here can be extended to more general matrix functions $f(A)$. Such an extension is likely not straightforward, as our bounds rely specifically on quasiseparable inversion formulas. Of course, in principle one could approximate $f(A)$ as a discretized Dunford-Cauchy integral in the form $f(A)=\sum_j w_i(z_jI-A)^{-1}$ and then apply decay bounds to each resolvent $(z_jI-A)^{-1}$, thus obtaining decay bounds for $f(A)$. On a deeper level, one might also think of combining a rational approximation for $f(z)$ with known results on the sum, product and inversion of quasiseparable matrices. On this topic, let us also recall that, under suitable hypotheses, functions of quasiseparable matrices retain an approximate quasiseparable structure: see for instance \cite{Gavrilyuk02, Massei17}.

\section*{Acknowledgments}
P. Boito acknowledges the MUR Excellence Department Project awarded to the Department of Mathematics, University of Pisa, CUP I57G22000700001. 

\noindent P. Boito is a member of GNCS-INdAM and is partially supported 
 by  European Union - NextGenerationEU under the National Recovery and Resilience Plan (PNRR) - Mission 4 Education and research - Component 2 From research to business - Investment 1.1 Notice Prin 2022 - DD N. 104  2/2/2022, titled Low-rank Structures and Numerical Methods in Matrix and Tensor Computations and their Application, proposal code 20227PCCKZ – CUP I53D23002280006.


\begin{thebibliography}{88}

\bibitem{Asplund59}
E.~Asplund, Inverses of matrices $\{a_{ij}\}$ which satisfy $a_{ij}= 0$ for $j> i+ p$. Mathematica Scandinavica (1959), 57--60.

\bibitem{Baskakov}
A.~G.~Baskakov. Estimates for the elements of inverse matrices, and the
spectral analysis of linear operators. Izv. Ross. Akad. Nauk Ser. Mat.,
61(6):3--26, 1997.

\bibitem{Bellavia}
S.~Bellavia, B.~Morini, and M.~Porcelli,  New updates of incomplete LU factorizations and applications to large nonlinear systems. Optimization Methods and Software, 29 (2014), 321--340.

\bibitem{BenziBoito14}
M.~Benzi and P.~Boito, Decay properties for functions of matrices over $C^*$-algebras. Linear Algebra and its Applications, 456 (2014), 174--198.

\bibitem{BenziRazouk07}
M.~Benzi and N.~Razouk, Decay bounds and $O(n)$ algorithms for approximating functions of sparse matrices. Electron. Trans. Numer. Anal 28 (2007), 16--39).

\bibitem{BenziRinelli}
M.~Benzi and M.~Rinelli, Refined decay bounds on the entries of spectral projectors associated with sparse Hermitian matrices. Linear Algebra and its Applications 647 (2022), 1--30.

%\bibitem{BenziRinelliSimunec}
%M.~Benzi, M.~Rinelli and I.~Simunec, Estimation of spectral gaps for sparse symmetric matrices (2024), %arXiv:2410.15349.

\bibitem{BenziSimoncini}
M.~Benzi and V.~Simoncini, Decay bounds for functions of Hermitian matrices with banded or Kronecker structure. SIAM Journal on Matrix Analysis and Applications, 36 (2015), 1263--1282.

\bibitem{BE23}
P.~Boito and Y.~Eidelman, Computation of quasiseparable representations of Green matrices. Linear Algebra Appl. 721 (2025), 47--80.

\bibitem{CH83}
C.~Chui and M.~Hasson, Degree of uniform approximation on disjoint intervals. Pacific Journal of Mathematics 105.2 (1983): 291--297.

\bibitem{Crouzeix}
M.~Crouzeix and C.~Palencia, The numerical range is a $(1+\sqrt{2})$-spectral set. SIAM Journal on Matrix Analysis and Applications, 38 (2017), 649--655.

\bibitem{SMC}
T.~A.~Davis and Y.~Hu, The University of Florida Sparse Matrix Collection. ACM Transactions on Mathematical Software 38, 1, Article 1 (December 2011), 25 pages. DOI: https://doi.org/10.1145/2049662.2049663

\bibitem{Dedieu88}
J.-P.~Dedieu, Matrix homographic iterations and bounds for the inverses of certain band matrices. Linear Algebra and its Applications 111 (1988): 29--42.

\bibitem{DMS84}
S.~Demko, W.~F.~Moss, and P.~W.~Smith,
Decay rates for inverses of band matrices. Mathematics of Computation 43.168 (1984), 491--499.

\bibitem{EGH1}
{Y.~Eidelman, I.~Gohberg and I.~Haimovici},
{\em Separable type representations of matrices and fast algorithms.
Volume 1. Basics. Completion problems. Multiplication and inversion 
algorithms}, Operator Theory: Advances and Applications,  Birkh\"auser, 2013.

\bibitem{FrommerNLAA}  
A.~Frommer, C.~Schimmel, and M.~Schweitzer, Bounds for the decay of the entries in inverses and Cauchy–Stieltjes functions of certain sparse, normal matrices. Numerical Linear Algebra with Applications 25.4 (2018): e2131.

\bibitem{Frommer18}
A.~Frommer, C.~Schimmel, M.~Schweitzer, Non-Toeplitz decay bounds for inverses of Hermitian positive definite tridiagonal matrices. Electronic Transactions on Numerical Analysis, 48 (2018), 362--372.

\bibitem{Gantmacher}
F.~R.~Gantmacher, 
{\em The Theory of Matrices. Vol. 1},
Chelsea Publishing Company, New York 1959.

\bibitem{Gavrilyuk02}
I.~P.~Gavrilyuk, W.~Hackbusch, and B.~N.~Khoromskij, $\mathcal{H}$-matrix approximation for the operator exponential with applications. Numerische Mathematik 92 (2002), 83--111.

\bibitem{Groechenig06}
K.~Gr\"ochenig and M.~Leinert, Symmetry and inverse-closedness of matrix algebras and functional calculus for infinite matrices. Transactions of the American Mathematical Society, 358 (2006), 2695--2711.

\bibitem{Hasson07}
M.~Hasson, The degree of approximation by polynomials on some disjoint intervals in the complex plane. Journal of Approximation Theory 144.1 (2007): 119--132.

\bibitem{Iserles}
A.~Iserles, How large is the exponential of a banded matrix?. University of Cambridge, Department of Applied Mathematics and Theoretical Physics (1999).

\bibitem{Jaffard}
S.~Jaffard, Propri\'et\'es des matrices ``bien localis\'es'' pr\`es de leur diagonale et quelques applications. In Annales de l'Institut Henri Poincar\'e C, Analyse non lin\'eaire (Vol. 7, No. 5, pp. 461-476), September 1990.

\bibitem{Krishtal}
I.~Krishtal, T.~Strohmer, and T.~Wertz, Localization of matrix factorizations. Foundations of Computational Mathematics, 15 (2015), 931--951.

\bibitem{Massei17}
S.~Massei and L.~Robol, Decay bounds for the numerical quasiseparable preservation in matrix functions. Linear Algebra and its Applications 516 (2017), 212--242.

\bibitem{Mastronardi}
N.~Mastronardi, M.~Ng, and E.~E.~Tyrtyshnikov, Decay in functions of multiband matrices. SIAM journal on matrix analysis and applications 31 (2010), 2721--2737.

\bibitem{Meurant}
G.~Meurant, A review on the inverse of symmetric tridiagonal and block tridiagonal matrices. SIAM Journal on Matrix Analysis and Applications 13 (1992), 707--728.

\bibitem{Nabben}
R.~Nabben, Decay rates of the inverse of nonsymmetric tridiagonal and band matrices. SIAM Journal on Matrix Analysis and Applications, 20 (1999), 820--837.

\bibitem{Pozza}
S.~Pozza and V.~Simoncini, Inexact Arnoldi residual estimates and decay properties for functions of non-Hermitian matrices. BIT Numerical Mathematics 59.4 (2019), 969--986.

\bibitem{Price51}
G.~B.~Price, Bounds for determinants with dominant principal diagonal. Proceedings of the American Mathematical Society, 2 (1951), 497--502.

%\bibitem{Schimmel}
%C.~Schimmel, {\em Bounds for the decay in matrix functions and its exploitation in matrix computations}, PhD Thesis, Wuppertal, Bergische Universit\"at, 2019.

\bibitem{Varah75}
J.~M.~Varah, A lower bound for the smallest singular value of a matrix. Linear Algebra Appl. 11 (1975), 3--5.

\end{thebibliography}
\end{document}